\documentclass[a4paper]{scrartcl} %
\usepackage{etex}

\PassOptionsToPackage{x11names}{xcolor}

\usepackage[utf8]{inputenc}
\usepackage[T1]{fontenc}
\usepackage[english]{babel}
\usepackage[neverdecrease]{paralist}
\usepackage{amssymb}
\usepackage{stmaryrd}
\usepackage[ntheorem]{empheq} %
\usepackage[thmmarks,amsmath]{ntheorem}
\usepackage[pdfusetitle,    colorlinks,
    linkcolor={red!70!black},
    citecolor={green!45!black},
    urlcolor={blue}]{hyperref}  %
\usepackage{jensntheorem}
\usepackage{lmodern}
\usepackage{microtype}
\usepackage{verbatim}
\usepackage{mathrsfs}
\usepackage{graphicx}
\usepackage{dsfont}
\usepackage{tikz-cd}
\usepackage{mathdots}
\usepackage[normalem]{ulem}

\numberwithin{equation}{section}

\usepackage{csquotes}
\usepackage[style=alphabetic,maxnames=10,maxalphanames=5,backend=biber]{biblatex}
\bibliography{prom}

\usetikzlibrary{arrows,positioning,decorations.markings}

\usepackage{jensstuff}

\title{Central leaves on Shimura varieties with parahoric reduction}
\author{Jens Hesse
\thanks{Technische Universit\"at Darmstadt, \texttt{math@jenshesse.eu}}}
\date{\today}

\begin{document}

\maketitle
	
\begin{abstract}
  We investigate the geometry of the special fiber of the integral model of a Shimura variety with parahoric level at a given prime place.

To be more precise, we deal with the definition of central leaves in this situation, their local closedness, and the relationship between the folations for varying parahoric level. This is connected to the verification of axioms for integral models formulated by He and Rapoport.
\end{abstract}

\tableofcontents

\addsec{Introduction}

Shimura varieties are objects of arithmetic geometry (namely varieties over number fields) that naturally arise in the search for generalized, non-abelian reciprocity laws (i.e., in the Langlands program) and as moduli spaces of abelian varieties (with certain extra structures on them). One way of approaching these objects is to try to understand their mod-$p$ reduction (which has to be carefully defined first). Insofar as a moduli interpretation in the above sense exists and continues to exist likewise for the mod-$p$ reduction\footnote{There need not be a \emph{literal} moduli interpretation, but in any event the stratifications in question derive from a close connection to moduli problems.}, it allows us to stratify the moduli space according to  several invariants of the abelian varieties parametrized, e.g., the isomorphism classes of their $p$-torsion. (An important observation is that these stratifications genuinely live in the characteristic $p$ world, making use of Frobenius endomorphisms and so on.) This, very roughly, is the general theme everything in this article revolves around.

More precisely, we will be dealing with Shimura varieties of Hodge type and parahoric level structure, at some fixed prime $v \mid p$ of the number field over which the Shimura variety is defined. Under some reasonably mild assumptions, cf.~\ref{std-assum}, Kisin and Pappas \cite{kisin-pappas} constructed a canonical integral model for such a Shimura variety. We try to understand some aspects of the geometry of the special fiber  of said integral model, namely the central leaves (roughly the patches where the isomorphism class of the $p$-divisible group associated with the abelian variety is constant); their being locally closed and how they vary as the (parahoric) level is varied.

Let us now go into more detail.

On the integral model $\sS_K$ ($K$ parahoric level) we have a ``universal'' abelian scheme (the quotation marks indicating that it is not really universal for some moduli problem on $\sS_K$, but it comes from a universal abelian scheme via pullback) and we have various kinds of Hodge tensors. We also have a ``universal'' isogeny chain of abelian schemes tightly connected to the ``universal'' abelian scheme. We define the \emph{naive central leaves} on the special fiber $\osS_K$ to be the loci where the isomorphism type of the geometric fibers of the $p$-divisible group associated with the ``universal'' abelian scheme (the ``universal'' $p$-divisible group) is constant (alternatively, we may use the ``universal'' isogeny chain). We arrive at the non-naive version by taking into account the Hodge tensors. This is the content of section~\ref{sec:defin-centr-leav}.

Next we show
\begin{TheoremA}\textnormal{(Corollary \ref{leaves-loc-closed})}
  The central leaves are locally closed.
\end{TheoremA}
using a somewhat simpler construction than the one given in \cite{hk}. By foundational work of Oort \cite{oort} we already know the naive central leaves to be locally closed. We show that the central leaves are open and closed inside the naive central leaves. Some basic topological considerations allow us to treat this question in perfected formal neighborhoods, allowing us to phrase it as a question about $p$-divisible groups with crystalline Tate tensors; a question which was answered by Hamacher \cite{hamacher}. This forms section~\ref{sec:local-clos-centr}.

Then we consider (tensor-respecting) self-quasi-isogenies of $p$-divisible groups (with tensors). The main take-away here is that, if we consider the geometric fibers of the ``universal'' isogeny chain of $p$-divisible groups in the above setting, then
\begin{TheoremA}\textnormal{(Example~\ref{QisgJbField})}
  The self-quasi-isogenies are independent of the level.
\end{TheoremA}
In section~\ref{sec:almost-product} we recall the almost product structure (under an additional technical assumption~\ref{assumpt-axiom-a} concerning the Rapoport-Zink uniformization map, which is satisfied e.g. if our reductive group over $\Q_p$ is residually split\footnote{A reductive group $G/\Q_p$ is \defn{residually split} if it has the same rank as its base change to the maximal unramified extension of $\Q_p$. The name ``residually split'' derives from the fact that in this case the maximal reductive quotients of parahoric group schemes associated with $G$ are split reductive.} \cite{zhou}), which expresses a variant of the Igusa variety (the \emph{Newton-Igusa variety}, with quasi-isogenies instead of isomorphisms) as a product of the Igusa variety and the Rapoport-Zink space.

We apply this to the question of how the central leaves behave under change of the parahoric level. We begin by constructing the change-of-parahoric morphisms on the Shimura varieties, Igusa varieties, Rapoport-Zink spaces and Newton-Igusa varieties. The almost product isomorphism is compatible with these morphisms (Lemma~\ref{cop-X}).

From this we can derive
\begin{TheoremA} \textnormal{(Corollary~\ref{cop-Ig-isom})}
  The map between Igusa varieties for varying parahoric level is an isomorphism.
\end{TheoremA}
This implies in particular
\begin{TheoremA}\textnormal{(Corollary~\ref{cop-surj})}
  The change-of-parahoric map between central leaves is surjective.
\end{TheoremA}
This is the most difficult part of the axioms on integral models He and Rapoport give in \cite{he-rapo}. Our proof is independent of the one by Rong Zhou in \cite{zhou}, even though we use some of the results he gave in the first version of the cited preprint, which did not contain a proof of the surjectivity. Positively answering another conjecture by He and Rapoport \cite[Rmk.~3.4]{he-rapo}, we show
\begin{TheoremA}\textnormal{(Corollary~\ref{cop-comp})}
  The change-of-parahoric map between central leaves is the composition of a flat universal homeomorphism of finite type and a finite Ã©tale morphism.
\end{TheoremA}

\subsection*{Acknowledgements}

This article essentially is an extract of my doctoral thesis \cite{diss} (another extract\footnote{In particular, there is a large overlap between the ``Background'' sections of the two articles.}, dealing with the EKOR stratification, is \cite{article-ekor}). I thank Torsten Wedhorn for suggesting the topic of the dissertation, his support and patiently answering my questions.
Moreover I thank Eva Viehmann and Paul Hamacher for their hospitality and helpful discussions during a month-long stay in Munich at the TU M\"unchen.
I am also grateful to Timo Richarz and Timo Henkel
for numerous helpful discussions.

This research was supported by the Deutsche Forschungsgemeinschaft (DFG), project number WE~2380/5.

\section{Background}
\label{cha:gener-prel}

\subsection{Shimura data of Hodge type}
\label{sec:shimura-data-hodge}

This article deals with aspects of the geometry of Shimura varieties (of Hodge type), which are the (systems of) varieties associated with Shimura data (of Hodge type).

\begin{Definition}\label{def-hodget}
  A \defn{Shimura datum of Hodge type}\index{Shimura datum of Hodge type} is a pair $(G,X)$, where $G$ is a reductive algebraic group over $\Q$ and $X\subseteq\Hom_{\R\text{-grp}}(\SSS,G_\R)$ is a $G(\R)$-conjugacy class ($\SSS:=\Res_{\C/\R}\G_{m,\C}$ being the Deligne torus\index{Deligne torus}) subject to the following conditions:
\begin{enumerate}[(1)]
\item For $h\in X$, the induced Hodge structure $\SSS\xrightarrow{h} G_\R \xrightarrow{\mathrm{Ad}}\GL(\Lie(G_\R))$ satisfies $\Lie(G_\C)=\Lie(G_\C)^{-1,1}\oplus\Lie(G_\C)^{0,0}\oplus\Lie(G_\C)^{1,-1}$.\label{sv1}
\item $\interior(h(i))\colon G^\mathrm{ad}_\R\to G^\mathrm{ad}_\R$ is a Cartan involution, i.e., $\{ g\in G^\mathrm{ad}(\C)\suchthat gh(i) = h(i)\overline{g}\}$ is compact. Another way of phrasing this condition: Every finite-dimensional real representation $V$ of $G^\mathrm{ad}_\R$ carries a $G^\mathrm{ad}_\R$-invariant bilinear form $\varphi$ such that $(u,v)\mapsto \varphi(u,h(i)v)$ is symmetric and positive definite. It is enough to show that this holds for one \emph{faithful} finite-dimensional real representation $V$.
\item $G^\mathrm{ad}$ \emph{cannot} be non-trivially written as $G^\mathrm{ad}\cong H\times I$ over $\Q$ with $\SSS\to G_\R \xrightarrow{\mathrm{proj}} H_\R$ trivial.
\item There exists an embedding $(G,X)\hookrightarrow (\GSp(V),S^\pm)$, where $(\GSp(V),S^\pm)$ is the Shimura datum associated with a finite-dimensional symplectic $\Q$-vector space $V$ (see below). That is, we have an embedding $G\hookrightarrow \GSp(V)$ of $\Q$-group schemes such that the induced map $\Hom_{\R\text{-grp}}(\SSS,G_\R)\hookrightarrow\Hom_{\R\text{-grp}}(\SSS,\GSp(V_\R))$ restricts to a map $X\hookrightarrow S^\pm$. \label{item:def-hodget-hodge-emb}
\end{enumerate}
\end{Definition}

\begin{Example}\label{sympl-ex}
  Let $W$ be a finite-dimensional $\R$-vector space.

  $\R$-group homomorphisms $\SSS\to \GL(W)$ then correspond to Hodge decompositions of $W$, i.e., to decompositions $W_\C=\oplus_{(p,q)\in\Z^2}W_\C^{p,q}$, such that $W_\C^{p,q}$ is the complex conjugate of $W_\C^{q,p}$ for all $(p,q)\in\Z^2$. Under this correspondence, $h\colon \SSS\to \GL(W)$ corresponds to the Hodge decomposition $W_\C^{p,q}=\{ w\in W_\C \suchthat \forall z\in\SSS(\R)=\C^\times\colon h(z)w=z^{-p}\bar{z}^{-q}w \}$. Hodge decompositions of $W$ of type $(-1,0)+(0,-1)$ correspond to complex structures on $W$: If $h\colon\SSS\to\GL(W)$ yields such a Hodge decomposition, then $h(i)$ gives an $\R$-endomorphism $J$ of $W$ with $J\circ J=-\id_W$.

  Let $V=(V,\psi)$ be a finite-dimensional symplectic $\Q$-vector space. We say that a complex structure $J$ on $V_\R$ is positive (resp. negative) if $\psi_J:=\psi_\R(\_,J\_)$ is a positive definite (resp. negative definite) symmetric bilinear form on $V_\R$. Define $S^+$ (resp. $S^-$) to be the set of positive (resp. negative) complex structures on $(V_\R,\psi_\R)$ and $S^\pm:=S^+\sqcup S^-$.

   We can make this more concrete: A symplectic basis of $(V_\R,\psi_\R)$ is a basis $e_1,\dotsc,\allowbreak e_g, e_{-g}, \dotsc, \allowbreak e_{-1}$, such that $\psi_\R$ is of the form
  \begin{equation*}
    \begin{pmatrix}
      & \tilde{I}_g \\
      -\tilde{I}_g &
    \end{pmatrix}
  \end{equation*}
  with respect to this basis, where $\tilde{I}_g=
  \begin{pmatrix}
    & & 1\\
    & \iddots & \\
    1 & &
  \end{pmatrix}$ is the antidiagonal identity matrix.\footnote{Occasionally (in particular when doing concrete matrix calculations), it is more convenient to number the basis vectors $1,\dotsc,g,-1,\dotsc,-g$ instead of $1,\dotsc,g,-g,\dotsc,-1$. Then the standard symplectic form is given by $\begin{smallpmatrix}
      & I_g \\
      -I_g &
    \end{smallpmatrix}$, $I_g$ being the $g\times g$ identity matrix.}

  Let $J$ be the endomorphism of $V_\R$ of the form
  \begin{equation*}
    \begin{pmatrix}
      & -\tilde{I}_g \\
      \tilde{I}_g &
    \end{pmatrix}
  \end{equation*}
  with respect to this basis. Then $J\in S^+$ and what we have described is a surjective map
  \begin{equation*}
    \{\text{symplectic bases of } (V_\R,\psi_\R)\} \twoheadrightarrow S^+.
  \end{equation*}

  In particular we see that $\Sp(V_\R,\psi_\R):=\{f\in\GL(V_\R) \suchthat \psi_\R(f(\_),f(\_))=\psi_\R\}$  (by virtue of acting simply transitively on the symplectic bases) acts transitively on $S^+\cong \Sp(V_\R,\psi_\R)/\SpO(V_\R,\psi_\R,J)$ (where we define $\SpO(V_\R,\psi_\R,J):=\Sp(V_\R,\psi_\R)\cap O(V_\R,\psi_J)=U((V_\R,J),\psi_J)$ for a fixed choice of $J\in S^+$) and therefore the general symplectic group $\GSp(V_\R,\psi_\R):=\{f\in\GL(V_\R) \suchthat \psi_\R(f(\_),f(\_))=c\cdot\psi_\R\text{ for some } c\in\R^\times\}$ acts transitively on $S^\pm$ (note that the element of the form $e_{\pm i}\mapsto e_{\mp i}$ of $\GSp(V_\R,\psi_\R)$ for any given choice of symplectic basis $\left(e_i\right)_i$ permutes $S^+$ and $S^-$).
\end{Example}

\begin{Definition}
  Condition~\eqref{sv1} of Definition~\ref{def-hodget} implies that the action of $\G_{m,\R}$ (embedded in $\SSS$ in the natural way) on $\Lie(G_\R)$ is trivial, so that $h$ induces a homomorphism ${w\colon \G_{m,\R}\to \Cent(G_\R)}$. This homomorphism is independent of the choice of $h\in X$ and is called the \defn{weight homomorphism} of $(G,X)$.

  Moreover, we denote by $\{\mu\}$ the the $G(\C)$-conjugacy class of the cocharacter $\mu_h:=h\circ(\id_{\G_{m,\C}},1)\colon \G_{m,\C}\to \G_{m,\C}^2\cong\SSS_\C\to G_\C$, where $h$ is as above. Obviously, the conjugacy class $\{\mu\}$ is independent of the particular choice of $h\in X$.
\end{Definition}

\begin{Remark}\label{conj-class-cochar}
  Let $L/\Q$ be a field extension such that $G_L$ contains a split maximal torus $T$.  Let $W:=\Norm_{G(L)}(T)/T$ be the Weyl group.  Then the natural map
  \begin{equation*}
    W\backslash \Hom_{L\text{-grp}}(\G_{m,L},T) \to G(L)\backslash \Hom_{L\text{-grp}}(\G_{m,L},G_L)
  \end{equation*}
  is bijective.

  Since the left hand side remains unchanged if we go from $L=\bar\Q$ (where as usual $\bar\Q$ denotes an algebraic closure of $\Q$) to $L=\C$, we see that $\{\mu\}$ contains a cocharacter defined over $\bar\Q$ and that we may then also consider $\{\mu\}$ as a $G(\bar\Q)$-conjugacy class.
\end{Remark}

\begin{Definition}
  The \defn{reflex field} $\bfE=\bfE(G,X)$ of $(G,X)$ is the field of definition of $\{\mu\}$, i.e., the fixed field in $\bar\Q$ of $\{\gamma\in\Gal(\bar\Q/\Q) \suchthat \gamma(\{\mu\})=\{\mu\}\}$.
\end{Definition}

\begin{Example}
  The reflex field of the Shimura datum $(\GSp_{2g,\Q},S^\pm)$ of Example~\ref{sympl-ex} is $\Q$. To wit, one of the cocharacters in the conjugacy class $\{\mu\}$ is
  \begin{equation*}
    \mu(z) =
    \begin{smallpmatrix}
      z & & & & & \\
    & \ddots & & & &\\
    & & z & & &\\
    & & & 1 & &  \\
    & & & & \ddots & \\
    & & & & & 1
    \end{smallpmatrix}.
  \end{equation*}
\end{Example}

\begin{Notation}
We denote the ring of (rational) adeles by $\A:=\A_\Q$, the subring of finite adeles by $\A_f:=\A_{\Q,f}$ and the subring of finite adeles away from some fixed prime $p$ by $\A_f^p$.
\end{Notation}

\begin{DefinitionRemark}\label{shvar}
  Let $K\subseteq G(\A_f)$ be a compact open subgroup. The \defn{Shimura variety of level $K$ associated with $(G,X)$} is the double coset space
  \begin{equation*}
    \Sh_K(G,X):=G(\Q)\backslash (X\times (G(\A_f)/K)).
  \end{equation*}

  A priori, this is just a set, but if $K$ is sufficiently small (i.e., ``neat'' in the sense of \cite{borelarith,pink}), $\Sh_K(G,X)$ can be canonically written as a finite disjoint union of hermitian symmetric domains.\footnote{If $K$ fails to be sufficiently small, one might very reasonably argue that our definition of the Shimura variety of level $K$ really is the definition of the \emph{coarse} Shimura variety and that one should be working with stacks instead.  Since we will only be interested in sufficiently small level, this is inconsequential for us.}  In particular, this gives $\Sh_K(G,X)$ the structure of a complex manifold.

  In fact, by the theorem of Baily-Borel, this complex manifold attains the structure of a quasi-projective complex variety in a canonical way.

  By work of Deligne, Milne and Borovoi, this variety is defined already (and again in a canonical way) over the reflex field $\bfE$. So in particular, it is defined over a number field independent of $K$.  This is important when varying $K$ and it is the reason why we consider the whole Shimura variety instead of its connected components over $\C$ on their own.  It is possible for the Shimura variety to have multiple connected components over $\C$ while being connected over $\bfE$.

  More detailed explanations may be found in \cite{milne-isv}.
\end{DefinitionRemark}

\subsection{Bruhat-Tits buildings}
\label{sec:bruh-tits-build}

Let $K$ be a complete discrete valuation field with ring of integers $\mathcal{O}$, uniformizer $\varpi$ and perfect residue field $\kappa:=\mathcal{O}/\varpi$.

\begin{Notation}
  For a (connected) reductive group $G$ over $K$, we denote by $\mathcal{B}(G,K)$ the extended (or enlarged) and by $\mathcal{B}^\mathrm{red}(G,K)$ the reduced (i.e., non-extended) Bruhat-Tits building of $G$ over $K$ \cite{bt-ii}. Moreover, $\mathcal{B}^\mathrm{abstract}(G,K)$ denotes the underlying abstract simplicial complex.
\end{Notation}

\begin{Remark}\label{gl-building}
  Let $V$ be a finite-dimensional $K$-vector space.
  
  As described in \cite[\nopp 1.1.9]{kisin-pappas} (originally in \cite{zbMATH03900941}), the points of $\mathcal{B}(\GL(V),K)$ correspond to graded periodic lattice chains $(\mathcal{L},c)$, i.e.,
\begin{itemize}
\item $\emptyset\ne\mathcal{L}$ is a totally ordered set of full $\mathcal{O}$-lattices in $V$ stable under scalar multiplication (i.e., $\Lambda\in\mathcal{L} \iff \varpi\Lambda\in\mathcal{L}$),
\item $c\colon \mathcal{L}\to\R$ is a strictly decreasing function such that $c(\varpi^n\Lambda)=c(\Lambda)+n$.
\end{itemize}
\end{Remark}

\begin{Remark}\label{gitterketten-nummerieren}
  Fix such an $\mathcal{L}$ and let $\Lambda^0\in\mathcal{L}$. Then every homothety class of lattices has a unique representative $\Lambda$ such that $\Lambda\subseteq\Lambda^0$ and $\Lambda\not\subseteq \varpi\Lambda^0$. Consider such representatives $\Lambda^i$ for all of the distinct homothety classes of lattices that make up $\mathcal{L}$. Because $\mathcal{L}$ is totally ordered and $\Lambda^i\not\subseteq \varpi\Lambda^0$, it follows that $\Lambda^i\supseteq \varpi\Lambda^0$ for all $i$ and that $\left\{\Lambda^i/\varpi\Lambda^0\right\}_i$ is a flag of non-trivial linear subspaces of $\Lambda^0/\varpi\Lambda^0\cong\kappa^{n}$, where $n:=\dim V$. Consequently, the number $r$ of homothety classes is in $\{1,\dotsc,n\}$; it is called the \defn{period length} (or \defn{rank}) of $\mathcal{L}$. Numbering the $\Lambda^i$ in descending order we hence obtain $r$ lattices $\Lambda^0,\Lambda^1,\dotsc,\Lambda^{r-1}$ such that
  \begin{equation}\label{eq:numbered-lattice-chain}
    \Lambda^0 \supsetneqq \Lambda^1 \supsetneqq \dotsb\supsetneqq \Lambda^{r-1} \supsetneqq \varpi\Lambda^0
  \end{equation}
  and $\mathcal{L}$ is given by the the strictly descending sequence of lattices
  \begin{equation*}
    \Lambda^{qr+i}=\varpi^q\Lambda^i,\quad q\in\Z, \; 0\leq i < r.
  \end{equation*}
\end{Remark}

\begin{Remark}\label{gsp-building}
  Let $V$ be a finite-dimensional symplectic $K$-vector space.
  
  $\mathcal{B}(\GSp(V),K)$ embeds into the subset of $\mathcal{B}(\GL(V),K)$ consisting of those $(\mathcal{L},c)$ such that
  $\Lambda\in\mathcal{L}\implies \Lambda^\vee\in \mathcal{L}$.

  Passing to the underlying abstract simplicial complexes means forgetting about the grading $c$ and
  \begin{equation*}
    \mathcal{B}^\mathrm{abstract}(\GSp(V),K) = \{\mathcal{L}\in\mathcal{B}^\mathrm{abstract}(\GL(V),K) \suchthat \Lambda\in\mathcal{L}\implies \Lambda^\vee\in \mathcal{L}\}.
  \end{equation*}

  If $\mathcal{L}\in\mathcal{B}^\mathrm{abstract}(\GSp(V),K)$ and $\{\Lambda^i\}_i$ is as in Remark~\ref{gitterketten-nummerieren}, then there is an involution $t\colon \Z\to\Z$ with $\left(\Lambda^i\right)^\vee=\Lambda^{t(i)}$, $t(i+qr)=t(i)-qr$, and $i< j\implies t(i)> t(j)$. So $-a:=t(0)> t(1)> \dotsb > t(r)=-a-r$, which implies $t(i)=-i-a$. Thus $i_0-t(i_0)=2i_0+a\in\{0,1\}$ for some unique $i_0\in\Z$. Hence, upon renumbering the $\Lambda^i$, we may assume that $a\in\{0,1\}$.

We therefore have
\begin{align*}
  \varpi\Lambda^0\subsetneqq\Lambda^{r-1}\subsetneqq \Lambda^{r-2}\subsetneqq \dotsb \subsetneqq \Lambda^0 \subseteq \left(\Lambda^{0}\right)^\vee=\Lambda^{-a} \subsetneqq \left(\Lambda^{1}\right)^\vee=\Lambda^{-1-a}  \\
  \subsetneqq \dotsb\subsetneqq \left(\Lambda^{r-1}\right)^\vee=\Lambda^{-r+1-a} \subseteq \Lambda^{-r}=\varpi^{-1}\Lambda^0.
\end{align*}
\end{Remark}

\subsection{Alteration of the Hodge embedding}
\label{sec:alt-hodge}

\begin{Notation}
  Let $E$ be a finite field extension of $\Q_p$.
  
  Denote by $\bE$ the completion of the maximal unramified extension of $E$ (hence $\bE=E\cdot\bQ_p$).
\end{Notation}

\begin{Remark}\label{ram-witt}
  If $E/\Q_p$ is unramified, then $\bOE=W(\bar\F_p)$, $\bar\F_p$ denoting an algebraic closure of $\F_p$ and $W\colon\mathrm{Ring}\to\mathrm{Ring}$ being the ($p$-adic) Witt vectors functor. This generalizes to the ramified case using \emph{ramified Witt vectors} instead, see e.g. \cite[Chap.~IV, (18.6.13)]{haze-fg} or \cite[Chapter~1]{ahsendorf}.
\end{Remark}

Let $(G,X)$ be a Shimura datum of Hodge type, let $(G,X)\hookrightarrow(\GSp(V),S^\pm)$ be an embedding as in Definition~\ref{def-hodget}\,\eqref{item:def-hodget-hodge-emb}, and let $x\in\mathcal{B}(G,\Q_p)$ be a point in the Bruhat-Tits building of $G$ over $\Q_p$.

We consider the associated Bruhat-Tits scheme ${\cal G}_x$, i.e., the affine smooth model of $G_{\Q_p}$ over $\Z_p$ such that ${\cal G}_x(\bZ_p)\subseteq G(\bQ_p)$ is the stabilizer of the facet of $x$ in ${\cal B}(G,\bQ_p)\overset{\text{\cite[Prop.~2.1.3]{landvogt}}}=\mathcal{B}(G,\uQ_p)$.
Let $K_p:={\cal G}_x(\Z_p)\subseteq G(\Q_p)$ and let $K^p\subseteq G(\A_f^p)$ be a sufficiently small open compact subgroup. Define $K:=K_pK^p\subseteq G(\A_f)$.

\begin{Assumptions}\label{std-assum}
  From now on, we will always make the following assumptions:
  \begin{itemize}
  \item $\mathcal{G}_x=\mathcal{G}_x^\circ$ is connected.
  \item $G$ splits over a tamely ramified extension of $\Q_p$.
  \item $p\nmid \#\pi_1(G^\mathrm{der})$.
  \end{itemize}
\end{Assumptions}

\begin{Notation}
  In order not to make notation overly cumbersome, we usually denote the base change $G_{\Q_p}$ of $G$ to $\Q_p$ by $G$ again.  (Later, we will almost exclusively be dealing with $G_{\Q_p}$.)
\end{Notation}

Under the above assumptions, Kisin and Pappas construct in \cite[section 1.2]{kisin-pappas} (building on \cite{landvogt}) a toral $G(\bQ_p)$- and $\Gal(\bQ_p/\Q_p)$-equivariant embedding
\begin{equation*}
  \iota\colon \mathcal{B}(G,\bQ_p) \to \mathcal{B}(\GL(V),\bQ_p),
\end{equation*}
restricting\footnote{In this case, this is obvious from Galois-equivariance, since $\mathcal{B}(G,K)=\mathcal{B}(G,K')^{\Gal(K'/K)}$ if $K'/K$ is an unramified extension.} to a map $\iota\colon \mathcal{B}(G,\Q_p)\to \mathcal{B}(\GL(V),\Q_p)$.

$\iota$ being a toral embedding means the following: $\iota$ is isometric after a suitable normalization of the norm on $\mathcal{B}(G,\bQ_p)$; moreover it is the canonical extension of a map ${\mathcal{B}^\mathrm{red}(G,\bQ_p)\to\mathcal{B}^\mathrm{red}(\GL(V),\bQ_p)}$ such that for each maximal $\bQ_p$-split torus $S\subseteq G$ there exists a maximal $\bQ_p$-split torus $T\subseteq \GL(V)$ such that $G \hookrightarrow \GL(V)$ maps $S$ into $T$ and $\iota$ restricts to a map between the reduced apartments associated with $(G,S)$ and $(\GL(V),T)$, respectively, compatible with translations.

This embedding $\iota$ depends on some choices:
By assumption, there exists a tamely ramified Galois extension $\tilde{K}/\Q_p$ with finite inertia group such that $G_{\tilde{K}}$ is split reductive. Let $H\to\Spec\Z_p$ be the split Chevalley form of $G$ over $\Z_p$. Then $\iota$ depends on the choice of
\begin{itemize}
\item an isomorphism $G_{\tilde{K}}\cong H_{\tilde{K}}$, 
\item a pinning $\left(T,M,f,R,\Delta,\left(X_\alpha\right)_{\alpha\in\Delta}\right)$ of $H$ (cf. \cite[Exp.~XXIII]{sga3}\footnote{Since $\Spec\Z_p$ is connected, some technicalities from \emph{loc. cit.} disappear here.}),
  which also entails the following:
Let $B$ be the Borel subgroup of $H$ corresponding to $\Delta$. When we talk about roots, it is with respect to $T$, and when we talk about positive roots and so on, it is with respect to $(T,B)$. We also fix a hyperspecial vertex $x_o$ in $\mathcal{B}(H,\Q_p)$ with stabilizer $H(\Z_p)$.
\item for every $\Q_p$-irreducible summand\footnote{Note that by reductivity and $\Char(\Q_p)=0$, every finite-dimensional representation of $G_{\Q_p}$ is completely reducible.} $V_i$ of the representation $G_{\Q_p}\to \GL(V)_{\Q_p}$ a lattice $\Lambda_i=U(\fn^-)v_i\subseteq V_i$, where $v_i\ne 0$ is a highest weight vector of $V_i$ and $\fn^-$ is the (strictly) negative root space inside the Lie algebra of $H$ over $\Z_p$,
\item a grading $c_{\Lambda_i}+t_i$ of the lattice chain $\{p^n\Lambda_i\}_{n\in\Z}$ given by real numbers $t_i\in\R$. Here ${(c_{\Lambda_i}+t_i)(p^n\Lambda_i)}:=n+t_i$.
\end{itemize}

By \cite[Lemma~2.3.3]{kisin-pappas}, we can (and do) arrange these choices in such a way that $\iota$ factors
\begin{equation*}
  \mathcal{B}(G,\bQ_p) \xrightarrow{j} \mathcal{B}(\GSp(V),\bQ_p) \to \mathcal{B}(\GL(V),\bQ_p),
\end{equation*}
where the last map is the canonical toral embedding (whose definition will be clear from Remark~\ref{gsp-building}). Again, $j$ restricts to a map
\begin{equation}\label{def-j}
  j\colon \mathcal{B}(G,\Q_p) \hookrightarrow \mathcal{B}(\GSp(V),\Q_p).
\end{equation}

Let $y=(\mathcal{L},c)$ be the image of $x\in\mathcal{B}(G,\Q_p)$ under the map \eqref{def-j} and let $\left(\Lambda^i\right)_i,r,a$ be as in Remark~\ref{gsp-building}. We define $N_p := \Stab_{\GSp(V)(\Q_p)}(\mathcal{L})$.

Consider the symplectic $\Q$-vector space
\begin{equation*} 
  V^\S := \bigoplus_{i=-(r-1)-a}^{r-1}V
\end{equation*}
(direct sum of symplectic spaces, i.e., if $\psi$ denotes the symplectic form on $V$, then the symplectic form $\psi^\S$ on $V^\S$ is given by $\bigoplus_{i=-(r-1)-a}^{r-1}\psi$) and the lattice in $V^\S_{\Q_p}$
\begin{equation*}
  \Lambda^\S:=\bigoplus_{i=-(r-1)-a}^{r-1}\Lambda^i.
\end{equation*}
By replacing $\Lambda^\S$ by a homothetic lattice, we may assume that $\Lambda^\S\subseteq\left(\Lambda^\S\right)^\vee$ (hence $\left(\Lambda^\S\right)^{\vee\vee}=\Lambda^\S$).

We have a diagonal embedding
\begin{equation*}
  \GSp(V) \hookrightarrow\GSp(V^\S)
\end{equation*}
and (with calligraphic script meaning that we talk about Bruhat-Tits group schemes)
\begin{equation*}
  \mathcal{GSP}(V)_y = \Stab_{\GSp(V)}(\mathcal{L}) = \bigcap\Stab_{\GSp(V)}(\Lambda^i) \subseteq \GSp(\Lambda^\S)\subseteq \GL(\Lambda^\S)\subseteq \mathcal{GL}(V)_{\Lambda^\S},
\end{equation*}
so that upon replacing our original embedding $G\hookrightarrow \GSp(V)$ by the embedding $G\hookrightarrow\GSp(V)\hookrightarrow\GSp(V^\S)$, we can assume that the parahoric subgroup on the Siegel side is given as
\begin{equation*}
  \Stab_{\GL(V^\S)(\Q_p)}(\Lambda^\S)\cap\GSp(V^\S)(\Q_p),
\end{equation*}
i.e., that our level is essentially given by a single lattice, albeit one that is (in general) not self-dual. This is \cite[\nopp 2.3.15]{kisin-pappas}.

\subsection{Siegel integral models}
\label{sec:siegel-integr-model}

With notation as above let
\begin{align*}
  N_p &:= \Stab_{\GSp(V)(\Q_p)}(\mathcal{L}) \quad\text{(as before)}, \\
  J_p&:= \Stab_{\GL(V^\S)(\Q_p)}(\Lambda^\S)\cap\GSp(V^\S)(\Q_p).
\end{align*}
Let $N^p\subseteq\GSp(V)(\A_f^p)$ and $J^p\subseteq\GSp(V^\S)(\A_f^p)$ be sufficiently small open compact subgroups, and $N:=N_pN^p$, $J:=J_pJ^p$.

In this subsection, we are going to describe integral models of $\Sh_{N}(\GSp(V),S^\pm)$ and of $\Sh_{J}(\GSp(V^\S),S^{\S,\pm})$ over $\Z_{(p)}$ and relate the two.

\begin{Remark}\label{rz-moduli}
  By \cite[Definition~6.9]{rz}, the integral model $\sS_{N}(\GSp(V),S^\pm)$ is given by the moduli problem $(\Z_{(p)}\text{-scheme})\ni S\mapsto \left\{(A,\bar\lambda,\eta^p)\right\}/\littlecong$, where:
  \begin{enumerate}[(a)]
  \item $A=\left(A_\Lambda\right)_{\Lambda\in\mathcal{L}}$ is an $\mathcal{L}$-set of abelian schemes, i.e.,
    \begin{itemize}
    \item for every $\Lambda\in\mathcal{L}$, an abelian $S$-scheme up to $\Z_{(p)}$-isogeny $A_\Lambda$ (i.e., $A_\Lambda$ is an object of the category $(\text{abelian } S\text{-schemes})\otimes\Z_{(p)}$, where the category $\mathcal{A}\otimes R$ for $\mathcal{A}$ an preadditive category and $R$ a ring has the same objects as $\mathcal{A}$ and $\Hom_{\mathcal{A}\otimes R}(X,Y)=\Hom(X,Y)\otimes_\Z R$ for all objects $X,Y$),
    \item for every inclusion $\Lambda_1\subseteq \Lambda_2$ a $\Z_{(p)}$-isogeny $\rho_{\Lambda_2,\Lambda_1}\colon  A_{\Lambda_1}\to A_{\Lambda_2}$,
    \item $\rho_{\Lambda_3,\Lambda_1}=\rho_{\Lambda_3,\Lambda_2}\circ\rho_{\Lambda_2,\Lambda_1}$ if $\Lambda_1\subseteq \Lambda_2 \subseteq \Lambda_3$ in $\mathcal{L}$,
    \item the height of $\rho_{\Lambda_2,\Lambda_1}$ is $\log_p|\Lambda_2/\Lambda_1|$. Here $\rho_{\Lambda_2,\Lambda_1}$ gives rise to a well-defined homomorphism of $p$-divisible groups, and what we mean is that the kernel of this homomorphism (which is a finite locally free commutative group scheme, which we also refer to simply as the kernel of $\rho_{\Lambda_2,\Lambda_1}$) is to have order $|\Lambda_2/\Lambda_1|$.
    \item For every $\Lambda\in \mathcal{L}$, there is an isomorphism (called \defn{periodicity isomorphism}) $\theta_\Lambda\colon A_\Lambda\to A_{p\Lambda}$ such that $\rho_{\Lambda,p\Lambda}\circ \theta_\Lambda = [p]\colon A_\Lambda\to A_\Lambda$ is the multiplication-by-$p$ isogeny.
    \end{itemize}
  \item $\bar\lambda\colon A\to\tilde{A}$ is a $\Q$-homogeneous principal polarization, i.e., a $\underline{\Q^\times}$-orbit of a principal polarization $\lambda\colon A\to \tilde{A}$. Here $\tilde{A}$ is the $\mathcal{L}$-set of abelian schemes over $S$ up to prime-to-$p$ isogeny given by $\tilde{A}_\Lambda:=(A_{\Lambda^\vee})^\vee$. And being a polarization $\lambda$ means being a quasi-isogeny of $\mathcal{L}$-sets $\lambda\colon A\to\tilde{A}$ such that
    \begin{equation*}
      A_\Lambda \xrightarrow{\lambda_\Lambda}\tilde{A}_\Lambda=(A_{\Lambda^\vee})^\vee\xrightarrow{\varrho_{\Lambda^\vee,\Lambda}^\vee}(A_\Lambda)^\vee
    \end{equation*}
    is a polarization of $A_\Lambda$ for all $\Lambda$. If $\lambda_\Lambda$ can be chosen to be an isomorphism up to prime-to-$p$ isogeny for all $\Lambda$, then we speak of a principal polarization. In that case, when referring to $\lambda_\Lambda$, we mean a $\lambda_\Lambda$ which is an isomorphism up to prime-to-$p$ isogeny.
  \item $\eta^p$ is a level-$N^p$-structure, i.e. (if $S$ is connected), it is a $\pi_1(S,s)$-invariant $N^p$-orbit of symplectic similitudes $\eta^p\colon V_{\A_f^p}\to H_1(A_s,\A_f^p)$ (where $s$ is some geometric basepoint and $H_1(A_s,\A_f^p)$ with its $\pi_1(S,s)$-action corresponds to the Tate $\A_f^p$-module of $A$ (cf. \cite[\nopp 6.8]{rz}), which is a smooth $\A_f^p$-sheaf). Note that this forces the abelian schemes $A_\Lambda$ to be $(\dim_\Q V)$-dimensional.
  \end{enumerate}
\end{Remark}

\begin{Definition}
  Set $\Lambda^\S_{\Z_{(p)}}:=\Lambda^\S_{\Z_p}\cap V^\S_\Q=\prod_{i=-(r-1)-a}^{r-1}\Lambda_{\Z_{(p)}}^i$. We choose a lattice $\Lambda^\S_\Z\subseteq V^\S$ such that $\Lambda^\S_\Z\otimes_\Z\Z_{(p)}=\Lambda^\S_{\Z_{(p)}}$ and $\Lambda^\S_\Z\subseteq (\Lambda^\S_\Z)^\vee$.
\end{Definition}

\begin{Remark}\label{paragraph-moduli}
  Set $d:=\bigl|\left(\Lambda_\Z^\S\right)^\vee/\Lambda_\Z^\S\bigr|$.  By \cite[\nopp 2.3.3, 3.2.4]{kisin}, the integral model $\sS_J(\GSp(V^\S),S^{\S,\pm})$ is given by the moduli problem $(\Z_{(p)}\text{-schemes})\ni S\mapsto \left\{(A^\S,\lambda^\S,\epsilon^p)\right\}/\littlecong$, where
  \begin{enumerate}[(a)]
  \item $A^\S$ is an abelian scheme over $S$ up to $\Z_{(p)}$-isogeny,
  \item $\lambda^\S\colon A^\S\to \left(A^\S\right)^\vee$ is a polarization of degree $d$ (i.e., the polarization of the (well-defined) associated $p$-divisible group has degree $d$),
  \item $\epsilon^p$ is a level-$J^p$-structure, i.e. (if $S$ is connected), it is a $\pi_1(S,s)$-invariant $J^p$-orbit of symplectic similitudes $\epsilon^p\colon V^\S_{\A_f^p}\to H_1(A^\S_s,\A_f^p)$. Note that this forces the abelian schemes $A^\S$ to be $(\dim_\Q V^\S)$-dimensional.
  \end{enumerate}
\end{Remark}

This completes the descriptions of the moduli problems, and we turn to the question of the relationship between the two.  Consider (for appropriate $N^p,J^p$; see below) the morphism $\chi\colon\sS_N(\GSp(V),S^\pm) \to \sS_J(\GSp(V^\S),S^{\S,\pm})$ given on $S$-valued points by sending $(A,\bar\lambda,\eta^p)$ to $(A^\S,\lambda^\S,\epsilon^p)$, where
\begin{enumerate}[(a)]
\item $\displaystyle A^\S:=\prod_{i=-(r-1)-a}^{r-1}A_{\Lambda^i}$,
\item $\displaystyle \lambda^\S:=\prod_{i=-(r-1)-a}^{r-1}\left(\rho_{\left(\Lambda^i\right)^\vee,\Lambda^i}^\vee\circ \lambda_{\Lambda^i}\right)$,
\item $\epsilon^p$ is the product $\prod_{i=-(r-1)-a}^{r-1}\eta^p$, to be interpreted as the product over $\eta^p\colon V_{\A_f^p}\to H_1(A_{\Lambda^i,s},\A_f^p)\cong H_1(A_s,\A_f^p)$, where the isomorphism $H_1(A_{\Lambda^i,s},\A_f^p)\cong H_1(A_s,\A_f^p)$ is by definition the identity for some fixed $i=i_0$ and otherwise induced by the transition map $\rho_{\Lambda^{i},\Lambda^{i_0}}$. We need that $N^p$ is mapped into $J^p$ by $\GSp(V)\hookrightarrow\GSp(V^\S)$ for this to make sense.
\end{enumerate}

\begin{Lemma}\label{tate-faithful}
  Let $S$ be a scheme, $\ell\ne p$ prime numbers. If $\ell$ does not appear as a residue characteristic of $S$, then the Tate module functors
  \begin{align*}
    H_1(\_,\Z_\ell)&\colon (\text{abelian } S\text{-schemes})\to (\text{Ã©tale }\Z_\ell\text{-local systems on } S), \\
    H_1(\_,\Q_\ell)&\colon (\text{abelian } S\text{-schemes})\to (\text{Ã©tale }\Q_\ell\text{-local systems on } S)
  \end{align*}
  (cf. \cite[\nopp III, 5.4 and 6.2]{groth-bt} for precise definitions) are faithful.

  If only $p$ and $0$ appear as residue characteristics of $S$, then the Tate module functor
  \begin{equation*}
    H_1(\_,\A_f^p)\colon (\text{abelian } S\text{-schemes})\to (\text{Ã©tale }\A_f^p\text{-local systems on } S)
  \end{equation*}
  is faithful.
\end{Lemma}

\begin{Proof}
  First note that the statements about $H_1(\_,\Q_\ell)$ and $H_1(\_,\A_f^p)$ follows from the statement about $H_1(\_,\Z_\ell)$, which is why it is enough to only look at $H_1(\_,\Z_\ell)$.

  A homomorphism of abelian $S$-schemes $f\colon A\to B$ vanishes if and only if it vanishes over every (geometric) fiber of $S$: Indeed, if it vanishes fiberwise, then it is flat by the fiber criterion for flatness. Applying that criterion again we see that the closed immersion and fiberwise isomorphism $\ker(f)\hookrightarrow A$ is flat, which means that is an isomorphism.

  This way we are reduced to the case where $R$ is an (algebraically closed) field of characteristic different from $\ell$. In this setting the faithfulness is well-known (the salient point being that the $\ell$-primary torsion is dense).
\end{Proof}

\begin{Lemma}\label{dantzig}
  Let $H$ be a totally disconnected locally compact\footnote{By (our) definition, locally compact implies Hausdorff.} group (i.e., a locally profinite group) and let $N\subseteq H$ be a compact subgroup. Then
  \begin{equation*}
    N = \bigcap_{\substack{N\subseteq J \\J\subseteq H\text{ open compact subgrp.}}} J.
  \end{equation*}
\end{Lemma}

Note that this is (a variant of) a well-known theorem by van Dantzig if $N=\{1\}$ \cite{vandantzig}.

\begin{Proof}
  We make use of the following fact \cite[Prop.~3.1.7]{arhangel}: A Hausdorff space is locally compact and totally disconnected if and only if the open compact sets form a basis of the topology. (Van Dantzig's theorem is the group version of this, which talks only about a neighborhood basis of the identity and open compact \emph{subgroups}.)

  First we show that $N$ is contained in some open compact subset $K\subseteq H$. For every $x\in N$ choose a compact open neighborhood $x\in K_x\subseteq H$. This is possible by the fact cited above. Then there is a finite subset $I\subseteq N$ such that $N\subseteq \bigcup_{x\in I}K_x=:K$.

  Next, for every $x\in N$ choose an open neighborhood of the identity $U_x$ such that $xU_xK\subseteq K$.  With $N\subseteq U:=\bigcup_{x\in N}xU_x$ we obtain $UK\subseteq K$. Replacing $U$ by $U\cap U^{-1}$, we may moreover assume it is symmetric. The subgroup generated by $U$ is open (hence closed) and contained in $K$, hence is an open compact subgroup.

  Thus $N$ even is contained in an open compact sub\emph{group}; in other words, we may assume that $H$ is compact, i.e., is a profinite group.

  Then $H/N$ is compact\footnote{Hausdorff quotient spaces of compact spaces are compact again, but for ``locally compact'' the analogous statement is not true in general!} and totally disconnected\footnote{Take $x,y\in H$ such that $xN\ne yN$. We show that any subspace $S\subseteq H/N$ containing both $xN$ and $yN$ is disconnected. Let $U\subseteq H/N$ be a neighborhood of $xN$ not containing $yN$. Let $x\in V\subseteq \pi^{-1}(U)$ be open and compact, where $\pi\colon H\to H/N$ is the projection. Then $yN\notin \pi(V)\subseteq H/N$ is open and compact (hence closed) and we have $S=(\pi(V)\cap S)\sqcup S\setminus \pi(V)$ where both $\pi(V)\cap S$ and $S\setminus\pi(V)$ are open in $S$. This shows that $S$ is disconnected.} (i.e., is a Stone space). By the fact cited above,
  \begin{equation*}
    H/N \supseteq\{1\} = \bigcap_{L\subseteq H/N\text{ open compact subset}} L.
  \end{equation*}

  Observe that the quotient map $H\to H/N$ is proper to deduce
  \begin{equation*}
    N = \bigcap_{\substack{N\subseteq M\\M\subseteq H\text{ open compact subset}}} M.
  \end{equation*}

  Say $M$ is an open and compact subset of $H$ containing $N$.  As we have shown above, there is an open compact subgroup $J\subseteq H$ in between $N$ and $M$, and this is all we need to complete the proof.
\end{Proof}

\begin{Proposition}\label{diag-emb}
  For every compact open subgroup $N^p\subseteq \GSp(V)(\A_f^p)$
  \begin{equation*}
    \chi\colon\sS_N(\GSp(V),S^\pm) \to \sS_J(\GSp(V^\S),S^{\S,\pm})
  \end{equation*}
  is a well-defined morphism for all compact open subgroups $N^p\subseteq J^p\subseteq \GSp(V^\S)(\A_f^p)$ and is a closed immersion for all sufficiently small compact open subgroups $N^p\subseteq J^p\subseteq \GSp(V^\S)(\A_f^p)$.
\end{Proposition}

\begin{Proof}
  The fact that it's well-defined is clear from the construction.

  To show the second statement, as in \cite[Prop.~1.15]{travaux}, it is enough to show that
  \begin{equation*}
    \sS_{N_pN^p}(\GSp(V),S^\pm) \to \varprojlim_{J^p}\sS_{J_pJ^p}(\GSp(V^\S),S^{\S,\pm})
  \end{equation*}
  is a closed immersion, i.e., a proper monomorphism.

  We begin by proving that it is a monomorphism, i.e., injective on $S$-valued points ($S$ arbitrary $\Z_{(p)}$-scheme). So, say $(A_1,\lambda_1,\eta_1^p)$ and $(A_2,\lambda_2,\eta_2^p)$ both map to $(A^\S,\lambda^\S,\epsilon_{J^p}^p)$. That means precisely that there is an isomorphism of abelian $S$-schemes up to $\Z_{(p)}$-isogeny
  \begin{equation*}
    \phi\colon \prod_{i=-(r-1)-a}^{r-1}A_{1,\Lambda^i} \xrightarrow{\cong} \prod_{i=-(r-1)-a}^{r-1}A_{2,\Lambda^i}
  \end{equation*}
  such that
  \begin{equation*}
    \phi^\vee \circ \prod_{i=-(r-1)-a}^{r-1}\left(\rho_{2,\left(\Lambda^i\right)^\vee,\Lambda^i}^\vee\circ \lambda_{2,\Lambda^i}\right) \circ \phi = \prod_{i=-(r-1)-a}^{r-1}\left(\rho_{1,\left(\Lambda^i\right)^\vee,\Lambda^i}^\vee\circ \lambda_{1,\Lambda^i}\right)
  \end{equation*}
  and
  \begin{equation*}
    H_1(\phi,\A_f^p)\circ \epsilon_{1,J^p}^p = \epsilon_{2,J^p}^p   \mod{J^p}.
  \end{equation*}

  We claim that $\phi$ comes from isomorphisms
  \begin{equation*}
    \phi_i\colon A_{1,\Lambda^i} \xrightarrow{\cong} A_{2,\Lambda^i}.
  \end{equation*}

  Certainly there is but one candidate for $\phi_i$: define $\phi_i$ to be the composition
  \begin{equation*}
    A_{1,\Lambda^i}\xrightarrow{\mathrm{incl}} \prod_{i=-(r-1)-a}^{r-1}A_{1,\Lambda^i} \xrightarrow{\phi} \prod_{i=-(r-1)-a}^{r-1}A_{2,\Lambda^i} \xrightarrow{\mathrm{proj}} A_{2,\Lambda^i}.
  \end{equation*}

  Our claim then is that
  \begin{equation*}
    \phi = \prod_{i=-(r-1)-a}^{r-1}\phi_i.
  \end{equation*}

  Apply $H^1(\_,\A_f^p)$ on both sides. For the left hand side, we have
  \begin{equation*}
    H_1(\phi,\A_f^p) = \epsilon_{2,J^p}^p\circ \left(\epsilon_{1,J^p}^p\right)^{-1}   \mod{J^p}.
  \end{equation*}
  and the right hand side of this equation is block diagonal. So
  \begin{equation*}
    H_1(\phi,\A_f^p) = \prod_{i=-(r-1)-a}^{r-1}H_1(\phi_i,\A_f^p) \mod{J^p}.
  \end{equation*}

  Since (by Lemma~\ref{dantzig})
  \begin{equation*}
    N^p = \bigcap_{\substack{N_\ell\subseteq J_\ell\\ J_\ell\subseteq \GSp(V^\S)(\Q_\ell)\text{ cpt. open subgrp.}}}J_\ell,
  \end{equation*}
  it follows that (with $\ell\ne p$)
  \begin{equation*}
    H_1(\phi,\Q_\ell) = \prod_{i=-(r-1)-a}^{r-1}H_1(\phi_i,\Q_\ell) \mod{N_\ell},
  \end{equation*}
  hence (since $N_\ell$ acts block-diagonally) that $H_1(\phi,\Q_\ell)= \prod_{i=-(r-1)-a}^{r-1}H_1(\phi_i,\Q_\ell)$.

  Since $H_1(\_,\Q_\ell)$ is faithful (Lemma~\ref{tate-faithful}), this implies $\phi=\prod_{i=-(r-1)-a}^{r-1}\phi_i$, as desired.

  Next, consider the extension by zero of $\left(H_1(\rho_{1/2,\Lambda^j,\Lambda^i},\A_f^p)\right)_{i,j}$ (where for ``$1/2$'' either ``$1$'' or ``$2$'' can be plugged in) to a map $H_1(A^\S,\A_f^p)\to H_1(A^\S,\A_f^p)$. Under the isomorphism given by the $J^p$-level structure this corresponds, up to the $J^p$-action, to the map $V^\S_{\A_f^p}\to V^\S_{\A_f^p}$ given by mapping the $i$'th copy of $V_{\A_f^p}$ identically to the $j$'th copy and the rest to zero. Thus $\rho_{1/2,i,j}$ yield the same up to $J^p$ after applying $H_1(\_,\A_f^p)$, hence they are equal in the $\Z_{(p)}$-isogeny category.

  Consequently, $\chi$ is a monomorphism.

  For properness, we will use the valuative criterion. Let $R$ be a discrete valuation ring with field of fractions $K$ and assume that a $K$-point $A^\S=\prod_{i=-(r-1)-a}^{r-1} A_{\Lambda^i}$ with its additional structures coming from $(A_{\Lambda^i})_i$ extends to an $R$-point $\mathcal{A}^\S$. Consider the map $A^\S\to A_{\Lambda^{i_0}}\to A^\S$ where the first map is a projection and the second an inclusion. By the NÃ©ron mapping property, this extends to a map $\mathcal{A}^\S\to\mathcal{A}^\S$. Define $\mathcal{A}_{\Lambda^{i_0}}$ to be the image of this map.

  The NÃ©ron mapping property also allows us to extend the transition isogenies $\rho_{\Lambda^{i_0},\Lambda^{j_0}}\colon\allowbreak {A_{\Lambda^{j_0}}\to A_{\Lambda^{i_0}}}$, $i_0\leq j_0$,  the periodicity isomorphisms, and the polarization.

  Since $\pi_1(\Spec K)$ surjects onto $\pi_1(\Spec R)$ (see \stacks{0BQM}), extending the level structure away from $p$ is trivial.
\end{Proof}

\subsection{Construction of the integral model}
\label{sec:constr-integr-model}

Let $\bfE$ be the reflex field of $(G,X)$ and $v\mid p$ a place of $\bfE$.

As the first step towards the construction of the integral model, we define $\sS_K^-(G,X)\to\Spec \mathcal{O}_{\bfE,(v)}$ to be the closure of $\Sh_K(G,X)$ in $\sS_N(\GSp(V),S^\pm)_{\mathcal{O}_{\bfE,(v)}}$ (or, equivalently by Proposition~\ref{diag-emb}, in $\sS_J(\GSp(V^\S),S^{\S,\pm})_{\mathcal{O}_{\bfE,(v)}}$). By this we mean the topological closure with the reduced subscheme structure or, equivalently (since $\Sh_K(G,X)$ is reduced), the flat closure (in the sense of \cite[\S\,2.8]{EGA4}).

Then we define $\sS_K(G,X)\to\Spec\mathcal{O}_{\bfE,(v)}$ to be the normalization of $\sS_K^-(G,X)\to\Spec\mathcal{O}_{\bfE,(v)}$.

This can be regarded as ``the obvious definition'', but it is entirely non-obvious and due to Kisin and Pappas that this really behaves as one would expect from a canonical integral model. These expectations, properly explained in \cite{kisin-pappas}, mainly concern the satisfaction of an extension property which is a weaker version of the valuative criterion for properness, and the existence of a local model diagram that essentially shows that the integral model is Ã©tale-locally isomorphic to the local model described by Pappas and Zhu \cite{pappas-zhu}. 

Define $E:=\bfE_v$, the $v$-adic completion of $\bfE$, and denote by $\kappa$ the residue field of $E$, a finite extension of $\F_p$. By abuse of notation, we also denote the base change of $\sS_K(G,X)$ to $\mathcal{O}_E$ by $\sS_K(G,X)$ again. By $\osS_K(G,X)$ we denote the base change of $\sS_K(G,X)$ to $\kappa$.

\subsection{Maps between Shimura varieties}
\label{sec:maps-between-shimura}

This section is based on \cite{zhou}, where more detailed explanations may be found. We let $K_p'=\mathcal{G}'(\Z_p)$ be another parahoric where $\mathcal{G}'$ is a Bruhat-Tits group scheme with $\mathcal{G}'=(\mathcal{G}')^\circ$. We assume that $K_p\subseteq K_p'$, i.e., that the facet (in the Bruhat-Tits building) associated with $K_p'$ is contained in the closure of the one associated with $K_p$.

\begin{Theorem} \textnormal{\cite[Theorem~7.1]{zhou}}\label{maps-sv}
  For sufficiently small $K^p$ there exists a morphism
  \begin{equation*}
    \pi_{K_pK^p,K_p'K^p}\colon \mathscr{S}_{K_pK^p}(G,X) \to \mathscr{S}_{K_p'K^p}(G,X).
  \end{equation*}
\end{Theorem}

In moving towards a proof, let us first note that $\GSp(V,\psi)\to\GSp(V^\S,\psi^\S)$ factors through
\begin{equation*}
  M := \left\{ \left(g_i\right)_i\in\prod_{i=-(r-1)-a}^{r-1}\GSp(V,\psi) \suchthat c(g_{-(r-1)-a})=\dotsb=c(g_{r-1}) \right\},
\end{equation*}
where $c\colon\GSp(V,\psi)\to \G_m$ is the multiplicator homomorphism.

There is a natural $X_M$ that makes $(M,X_M)$ into a Shimura datum. Define
\begin{equation*}
  H_p:=\Stab_{M(\Q_p)}(\Lambda^\S), \quad J_p:=\Stab_{\GSp(V^\S_{\Q_p},\psi^\S_{\Q_p})}(\Lambda^\S).
\end{equation*}

Then for sufficiently small $H^p,J^p$, we obtain a closed immersion
\begin{equation*}
  i\colon \Sh_{H_pH^p}(M,X_M) \hookrightarrow \Sh_{J_pJ^p}(\GSp(V^\S,\psi^\S),S^{\S,\pm}).
\end{equation*}

$\Sh_{H_pH^p}$ has a moduli interpretation over $\Z_{(p)}$ (essentially: $\#\{-(r-1)-a,\dotsc,r-1\}=2(r-1)+a+1$ abelian schemes up to prime-to-$p$ isogeny endowed with certain polarizations and prime-to-$p$ level structure).

\begin{Proposition} \textnormal{\cite[Prop.~7.2]{zhou}}\label{maps-sv-climm}
  If $J^p$ is sufficiently small, then $i$ extends to a closed immersion
  \begin{equation*}
    i \colon \mathscr{S}_{H_pH^p}(M,X_M)
    \to \mathscr{S}_{J_pJ^p}(\GSp(V^\S,\psi^\S),S^{\S,\pm}).
  \end{equation*}
\end{Proposition}

Now consider the embedding $i\colon \mathcal{B}(G,\Q_p) \to \mathcal{B}(\GSp(V,\psi),\Q_p)$.

We still have $K_p=\mathcal{G}(\Z_p)$ with $\mathcal{G}=\mathcal{G}_x$, $x\in\mathcal{B}(G,\Q_p)$. Let $\fg$ be the minimal facet in $\mathcal{B}(\GSp(V,\psi),\Q_p)$ containing $i(x)$. So $\fg$ corresponds to some lattice chain $\Lambda^0\supseteq\dotsb\supseteq\Lambda^{r-1}$ as above. We have
\begin{equation*}
  \Sh_{K_pK^p}(G,X)\xrightarrow{\text{cl. imm.}}\Sh_{H_pH^p}(M,X_M)_\bfE
  \xrightarrow{\text{cl. imm.}}\Sh_{J_pJ^p}(\GSp(V^\S,\psi^\S),S^{\S,\pm})_\bfE
\end{equation*}
and $\mathscr{S}^-_K(G,X)$ is defined to be the closure of $\Sh_{K_pK^p}(G,X)$ in $\mathscr{S}_{J_pJ^p}(\GSp(V^\S,\psi^\S),S^{\S,\pm})_{\mathcal{O}_{\bfE,(v)}}$.

\begin{Corollary} \textnormal{\cite[Cor.~7.3]{zhou}}\label{maps-sv-closure}
  $\mathscr{S}^-_K(G,X)$ also can be described as the closure of $\Sh_{K_pK^p}(G,X)$ in \linebreak[1] $\mathscr{S}_{H_pH^p}(M,X_M)_{\mathcal{O}_{\bfE,(v)}}$.
\end{Corollary}

Let $\ff'$ be a facet of $\mathcal{B}(\GSp(V^\S,\psi^\S),\Q_p)$ in the closure of $\fg$. Then $\ff'$ corresponds to a ``sub-lattice chain'' $\Lambda^{i_1}\supseteq\dotsb\supseteq\Lambda^{i_s}$, $\left\{i_1,\dotsc,i_s\right\}\subseteq\{0,\dotsc,r-1\}$. Defining $M',H_p'$ as above, we naturally obtain morphisms $M\to M'$ and 
\begin{equation*}
  \omega_{H_pH^p,H_p'\left.H'\right.^p}\colon \mathscr{S}_{H_pH^p}(M,X_M)\to \mathscr{S}_{H'_p\left.H'\right.^p}(M',X_{M'})
\end{equation*}
for suitable levels $H^p,\left.H'\right.^p$ away from $p$.

\begin{ProofSketch}[of Theorem~\ref{maps-sv}]
  Let $\ff$ be the facet of $K_p$ (that is to say $\ff$ is the minimal facet satisfying $x\in \ff$), $\ff'$ that of $K_p'$. Then $\ff\subseteq\overline{\ff'}$.

  Let $x'\in\ff'$ be so close to $x$ that $\fg\subseteq\overline{\fg'}$ in $\mathcal{B}(\GSp(V,\psi),\Q_p)$, where $\fg$ and $\fg'$ denote the minimal facets containing $i(x)$ and $i(x')$ respectively.

  The constructions from above yield:
  \begin{equation}\label{eq:incomplete-map-diagram}
    \begin{tikzpicture}[baseline=(current  bounding  box.center),node distance=4cm, auto]
      \node (X) at (0,3) {$\mathscr{S}_{K_pK^p}^-(G,X)$};
      \node (Z) at (7,3) {$\mathscr{S}_{K'_p\left.K'\right.^p}^-(G,X)$};
      \node (Zs) at (7,0) {$\mathscr{S}_{H_p'\left.H'\right.^p}(M',X_{M'})_{\mathcal{O}_{\bfE,(v)}}$};
      \node (Xs) at (0,0) {$\mathscr{S}_{H_p\left.H\right.^p}(M,X_{M})_{\mathcal{O}_{\bfE,(v)}}$};
      
      \draw[right hook->] (X) to node[swap] {closed imm.} (Xs);
      \draw[right hook->] (Z) to node {closed imm.} (Zs);
      \draw[->] (Xs) to node {$\omega_{H_pH^p,H_p'\left.H'\right.^p}$} (Zs);
    \end{tikzpicture}
  \end{equation}

  On the generic fiber we can complete this to a commutative square.

  By Corollary~\ref{maps-sv-closure}, this implies that \eqref{eq:incomplete-map-diagram} also completes to a commutative square. Now normalize.
\end{ProofSketch}

\subsection{Local structure of the integral model}
\label{sec:local-structure}

\subsubsection{Generizations and irreducible components}
\label{sec:gener-irred-comp}

Let $\mathscr{X}\to \Spec \mathcal{O}_{\bE}$ be a flat
scheme locally of finite type; denote the special fiber by $X\to \Spec\bar\F_p$ and the generic fiber by $\mathcal{X}\to\Spec\bE$. We assume that $\mathcal{X}$ is locally integral (e.g. smooth).

For example, we can consider $(\mathscr{X},X,\mathcal{X})=(\mathscr{S}^-_K(G,X)_{\bOE}, \mathscr{S}^-_K(G,X)_{\bOE}\otimes_{\bOE}\bar\F_p,\allowbreak {\Sh_K(G,X)\otimes_E\bE})$.

Let $\bar x\in X(\bar\F_p)$. 

\begin{Lemma}\label{closed-generizations}
  There is a generization $x$ of $\bar x$ which lies in the generic fiber $\mathcal{X}$, and is a closed point in there, i.e., $x\in \mathcal{X}(L)$ for a finite extension $L/\bE$.
\end{Lemma}

\begin{Definition}\label{def-cpg}
  We shall call such a point $x$ a \defn{closed point generization} of $\bar x$ for short.
\end{Definition}

\begin{Proof}
  Due to flatness (going-down) there is \emph{some} generization in the generic fiber; call it $x_0$.

  By \stacks{053U} the following set is dense (and in particular non-empty) in the closure of $\{x_0\}$ in $\mathcal{X}$:
  \begin{equation*}
    \left\{ x\in\mathscr{X} \suchthat x\text{ is a specialization of } x_0\text{ and a closed point generization of } \bar x \right\}.
  \end{equation*}
\end{Proof}

\begin{Lemma}\label{closed-generizations-realization}
  Notation as in the preceding lemma.

  The specialization $x\leadsto \bar x$ can be realized by an ${\cal O}_L$-valued point of $\mathscr{X}$.
\end{Lemma}

\begin{Proof}
  First off, by \cite[\nopp 7.1.9]{EGA2}, it can be realized by a morphism $\Spec R=\{\eta,s\}\to\mathscr{X}$ of $\bOE$-schemes, where $R$ is a discrete valuation ring such that $L\cong\kappa(\eta)=\Quot(R)$ as field extensions of $\kappa(x)$.

  We hence get local homomorphisms of local rings $\bOE\to{\cal O}_{\mathscr{X},\bar x}\to R$.

  Thus the discrete valuation on $L$ defined by $R$ extends the discrete valuation on $\bE$. But there is but one such extension and its valuation ring is ${\cal O}_L$ (by definition).
\end{Proof}

\begin{Lemma}
  Mapping $x$ to the unique irreducible component of $\mathscr{X}$ that contains $x$ establishes a surjection from the set of closed point generizations $x$ of $\bar x$ to the set of irreducible components of $\mathscr{X}$ containing $\bar x$.
\end{Lemma}

\begin{Proof}
  If $x_0\in\mathcal{X}$ is a generization of $\bar x$, then $x_0$ lies in a unique irreducible component of $\mathscr{X}$ because $\mathcal{X}$ is locally irreducible. Hence the map described above is well-defined.

  Now for surjectivity: Given an irreducible component $C$ of $\mathscr{X}$ containing $\bar x$, let $x_0\in C$ be the generic point. Then $x_0$ must be in the generic fiber (else we would be able to find a generization in the generic fiber by going-down). Now go through the proof of Lemma~\ref{closed-generizations} with this particular choice of $x_0$.
\end{Proof}

\subsubsection{Normalization and completion}
\label{sec:norm-compl}

For reference, we collect some facts concerning the passage to normalization and completion and in particular how it applies to integral models of Shimura varieties.

\begin{Fact}
  $\mathscr{S}^-_K(G,X) \to \Spec\bOE$ is quasi-projective, so in particular of finite type.

  Hence $\mathscr{S}^-_K(G,X)$ and $\mathscr{S}^-_K(G,X)_{\bOE}$ are excellent.

  As a consequence the normalization
  \begin{equation*}
    \mathscr{S}_K(G,X)_{\bOE}\xrightarrow{\nu}\mathscr{S}^-_K(G,X)_{\bOE}
  \end{equation*}
  is finite. (This really is a normalization because normalization and completion behave well together in the excellent case (to get from ${\cal O}_{E^\mathrm{ur}}$ to $\bOE$ and from ${\cal O}_{\bfE,(v)}$ to ${\cal O}_{E}$) and because normalization commutes with base change along filtered colimits of smooth morphisms (to get from ${\cal O}_{E}$ to ${\cal O}_{E^\mathrm{ur}}$)).
\end{Fact}

We will always denote by $(\;)^\sim$ the integral closure of a ring in its total ring of fractions.

$\bar x$ still denotes an $\bar\F_p$-valued point of $\mathscr{S}^-_K(G,X)$. Let $\nu^{-1}(\bar x)=\{\bar y=\bar{y}_1,\dotsc,\bar{y}_n\}$.

Let $\mathscr{S}^-:=\mathscr{S}^-_K(G,X)_{\bOE}$ and $\mathscr{S}:=\mathscr{S}_K(G,X)_{\bOE}$.

\begin{Fact}
  $\left(\nu_*{\cal O}_{\mathscr{S}^-_K(G,X)}\right)_{\bar x}^{\wedge}=\prod_{j=1}^n \widehat{\cal O}_{\mathscr{S},\bar{y}_j}$.
\end{Fact}

\begin{Fact} By \stacks{0C3B}:
  $\left(\nu_*{\cal O}_{\mathscr{S}^-_K(G,X)}\right)_{\bar x}=\widetilde{\cal O}_{\mathscr{S}^-,\bar x}$.
\end{Fact}

\begin{Fact} By \stacks{035P}:
  $\widetilde{A} = \prod_{\fq\in\Min(A)}(A/\fq)^{\sim}$, if $A$ is a reduced ring and $\#\Min(A)<\infty$.
\end{Fact}

\begin{Fact}
  $(\;)^\sim$ and $(\;)^\wedge$ commute in the case of excellence (see e.g. \cite[7.6.1, 7.8.3.1(vii)]{EGA4}).

  For instance, $\left.{\cal O}_{\mathscr{S}^-,\bar x}^\wedge\right.^\sim=\left.{\cal O}_{\mathscr{S}^-,\bar x}^\sim\right.^\wedge$ (on the right hand side, we have a completion of a ${\cal O}_{\mathscr{S}^-,\bar x}$-module).
\end{Fact}

\begin{Fact}
  In the case of excellence, the completion of a normal domain is a normal domain (see e.g. \cite[\nopp 12.50]{gw}).
\end{Fact}

Thus we have
\begin{equation}\label{eq:1}
  \begin{aligned}
    \prod_{\fq\in\Min(\widehat{\cal O}_{\mathscr{S}^-,\bar x})}(\widehat{\cal O}_{\mathscr{S}^-,\bar x}/\fq)^\sim
    &\cong \left.\widehat{\cal O}_{\mathscr{S}^-,\bar x}^\wedge\right.^\sim
    \cong \left.\widehat{\cal O}_{\mathscr{S}^-,\bar x}^\sim\right.^\wedge \\
    &\cong \left(\nu_*{\cal O}_{\mathscr{S}^-_K(G,X)}\right)_{\bar x}^{\wedge}
    \cong \prod_{j=1}^n \widehat{\cal O}_{\mathscr{S},\bar{y}_j}
  \end{aligned}
\end{equation}
and the rings $\widehat{\cal O}_{\mathscr{S},\bar{y}_j}$ are normal domains. Hence we obtain a bijection
\begin{equation*}
  \Min(\widehat{\cal O}_{\mathscr{S}^-,\bar x}) \xleftrightarrow{\;1:1\;}\nu^{-1}(\bar x)
\end{equation*}
such that there exists a numbering $\fq_1,\dotsc,\fq_n$ of the elements of $\Min(\widehat{\cal O}_{\mathscr{S}^-,\bar x})$ such that \eqref{eq:1} restricts to an isomorphism
\begin{equation*}
  \widehat{\cal O}_{\mathscr{S}^-,\bar x}/\fq_j \cong \widehat{\cal O}_{\mathscr{S},\bar{y}_j}
\end{equation*}
(also see \cite[\nopp 7.6.2]{EGA4}).

Also:
\begin{align*}
  \left.\widehat{\cal O}_{\mathscr{S}^-,\bar x}^\sim\right.^\wedge
  &\cong \left(\prod_{\fq\in\Min({\cal O}_{\mathscr{S}^-,\bar x})}({\cal O}_{\mathscr{S}^-,\bar x}/\fq)^\sim\right)^\wedge \\
  &\cong \prod_{\fq\in\Min({\cal O}_{\mathscr{S}^-,\bar x})}\left.\left({\cal O}_{\mathscr{S}^-,\bar x}/\fq\right)^\sim\right.^\wedge \\
  &\cong \prod_{\fq\in\Min({\cal O}_{\mathscr{S}^-,\bar x})}\left.\left({\cal O}_{\mathscr{S}^-,\bar x}/\fq\right)^\wedge\right.^\sim \\
  &\overset{\mathclap{\text{Artin-Rees}}}\cong \prod_{\fq\in\Min({\cal O}_{\mathscr{S}^-,\bar x})}\left(\widehat{\cal O}_{\mathscr{S}^-,\bar x}/\widehat\fq\right)^\sim, \quad \widehat\fq=\fq\widehat{\cal O}_{\mathscr{S}^-,\bar x}.
\end{align*}

Now by \cite[\nopp 2.1.2, 4.2.2]{kisin-pappas} for all $\fq\in\Min(\widehat{\cal O}_{\mathscr{S}^-,\bar x})$ we have that
\begin{equation*}
  \left(\widehat{\cal O}_{\mathscr{S}^-,\bar x}/\fq\right)^\sim=\widehat{\cal O}_{\mathscr{S}^-,\bar x}/\fq\cong R_{G,\bar x,\fq}=R_G
\end{equation*}
is normal and $\widehat{\cal O}_{\mathscr{S}^-,\bar x}/\fq\cong \widehat{\cal O}_{\mathscr{S},\bar y}$ for an appropriate choice of $\fq$ (i.e., of $x$). (The notation $R_G$ is from \cite{kisin-pappas}, where it is defined as the formal local ring of the local model.)

\subsection{Hodge tensors and (lack of) moduli interpretation}

We discuss the moduli interpretation of Shimura varieties of Hodge type and the partial extension of this interpretation to the integral model as constructed in Section~\ref{sec:constr-integr-model}.
Let $(G,X)$ be a Shimura datum of Hodge type, let $(G,X)\hookrightarrow(\GSp(V),S^\pm)$ be an embedding as in Definition~\ref{def-hodget}\,\eqref{item:def-hodget-hodge-emb}.

\subsubsection{The story for the \texorpdfstring{$\C$}{C}-valued points}
\label{sec:story-c-valued}

\begin{Lemma} \textnormal{\cite[Prop.~3.1]{deligne-hc}}\label{hodge-tensors-rat}
  There exist numbers $n,r_i,s_i\in\Z_{\geq 0}$ and tensors $t_i\in V^{\otimes r_i}\otimes (V^*)^{\otimes s_i}$, $1\leq i\leq n$, such that $G$ is the subgroup of $\GSp(V)$ fixing the $t_i$.
\end{Lemma}

\begin{Remark} (See \cite[Section~7]{milne-isv}.)\label{hodge-moduli-linear}
  For $K$ a compact open subgroup of $G(\A_f)$, the set $\Sh_K(G,X)(\C)$ of Definition~\ref{shvar} has the following moduli interpretation: isomorphism classes of triples $((W,h),\left(u_i\right)_{0\leq i\leq n},\eta K)$, where
  \begin{enumerate}[(a)]
  \item $(W,h)$ rational Hodge structure of type $(-1,0)+(0,-1)$,
  \item $\pm u_0$ polarization of $(W,h)$ (i.e., either $u_0$ or $-u_0$ is a polarization),
  \item $u_i\in V^{\otimes r_i}\otimes (V^*)^{\otimes s_i}$ for $1\leq i\leq n$,
  \item $\eta K$ is a $K$-orbit of isomorphisms $V\otimes \A_f\xrightarrow{\sim} W\otimes\A_f$, mapping $\psi$ to a $\A_f^\times$-multiple of $u_0$, and every $t_i$ to $u_i$ ($i\geq 1$),
  \end{enumerate}
  such that there exists an isomorphism $a\colon W \xrightarrow{\sim} V$ mapping $u_0$ to a $\Q^\times$-multiple of $\psi$, every $u_i$ to $t_i$ ($i\geq 1$), and $h$ to an element of $X$ (so $^ah:=\GL(a)\circ h\colon\SSS\to\GL(V)$ is in $X$).
\end{Remark}

\begin{ProofSketch}
  Given $((W,h),\left(u_i\right)_{0\leq i\leq n},\eta K)$ and $a$ as above, we consider the pair ${(^a,a\circ\eta)}$. By assumption, $^ah\in X$ and $a\circ \eta$ is a symplectic similitude fixing the $t_i$; hence $(^ah,a\circ \eta)\in X\times G(\A_f)$. The double quotient now results precisely from the ambiguity in the choices of $a$ and $\eta\in\eta K$.
\end{ProofSketch}

\begin{Remark}
  Denote by $\Q(r)$ the rational Hodge structure of type $(-r,-r)$ with underlying vector space $\Q$, and
  denote the multiplier character $\GSp(V)\to\G_m$ as well as its restriction to $G$ by $c$. Then for $h\in X$, $c\circ h\cong \Q(1)$ as rational Hodge structure, and the symplectic form gives an isomorphism $V\cong V^*\otimes\Q(1)$.

  With this, we may also interpret a tensor $t\in V^{\otimes r}\otimes (V^*)^{\otimes s}$ as a multilinear map $V^{r+s}\to\Q$, and, if $t$ is fixed by $G$, as a morphism $(V,h)^{\otimes(r+s)}\to\Q(r)$ of Hodge structures. Provided $t\ne 0$, this implies $r+s=2r$, i.e. $r=s$.
\end{Remark}

\begin{Remark}
  Since $A\mapsto H_1(A(\C),\Q)$ yields an equivalence between the category of complex abelian varieties up to isogeny and the category of polarizable rational Hodge structures of type $(-1,0)+(0,-1)$, we may also take $\Sh_K(G,X)(\C)$ to be a moduli problem of abelian varieties with extra structure: consider isomorphism classes of triples $(A,\left(u_i\right)_{0\leq i\leq n},\eta K)$, where
  \begin{enumerate}[(a)]
  \item $A$ is a complex abelian variety up to isogeny,
  \item $\pm u_0$ is a polarization of the rational Hodge structure $V:=H_1(A(\C),\Q)$,
  \item $u_1,\dotsc,u_n$ is as in Remark~\ref{hodge-moduli-linear},
  \item $\eta K$ is a $K$-orbit of $\A_f$-linear isomorphisms $V\otimes \A_f\to V_f(A):=T_f(A)\otimes_\Z\Q=\left(\varprojlim A(k)[n]\right)\otimes_\Z\Q$, mapping $\psi$ to a $\A_f^\times$-multiple of $u_0$ and every $t_i$ to $u_i$ ($i\geq 1$),
  \end{enumerate}
  such that there exists an isomorphism $a\colon H_1(A(\C),\Q)\xrightarrow{\sim} V$ mapping $u_0$ to a $\Q^\times$-multiple of $\psi$, $u_i$ to $t_i$ ($i\geq 1$) and $h$ to an element of $X$.
\end{Remark}

\begin{Remark}\label{univ-ab-sch-rat}
  When rephrasing the moduli problem appropriately, the Shimura variety as an $\bfE$-variety is given by such a moduli problem, see \cite[Section~14]{milne-isv}.

  In particular, we have a universal abelian scheme $\mathcal{A}\to\Sh_K(G,X)$.
\end{Remark}

\subsubsection{The story for the integral model}
\label{sec:story-integral-model}

Integral models of Hodge type Shimura varieties in general seem not to afford nice straightforward moduli interpretations. Still, the moduli interpretation of the $\C$-valued points extends to \emph{some} degree. First, we do have a generalization of Lemma~\ref{hodge-tensors-rat} as follows.

\begin{Proposition}\label{hodge-tensors-kisin} \textnormal{\cite[Prop.~1.3.2]{kisin}} Let $R$ be a discrete valuation ring of mixed characteristic, $M$ a finite free $R$-module, and $\mathcal{G}\subseteq\GL(M)$ a closed $R$-flat subgroup with reductive generic fiber. Then $\mathcal{G}$ is defined by a finite collection of tensors $\left(s_\alpha\right)_\alpha\subset M^\otimes$.
\end{Proposition}

Here $M^\otimes$ is the direct sum of all the $R$-modules which can be formed from $M$ by taking duals and (finite) tensor products.\footnote{Kisin in \emph{loc. cit.} also allows for taking symmetric and exterior powers. As Deligne pointed out, this is unnecessary. \cite{deligne-letter-kisin}}

We return to the setting of Section~\ref{sec:alt-hodge}; in particular we consider a Bruhat-Tits scheme $\mathcal{G}:=\mathcal{G}_x$.

By the proposition just stated, we find a collection of tensors $\left(s_\alpha\right)_\alpha\subset (\Lambda^\S_{\Z_{(p)}})^\otimes$ whose pointwise stabilizer is the Zariski closure $G_{\Z_{(p)}}$ of $G$ in $\GL(\Lambda^\S_{\Z_{(p)}})$.

\begin{Remark}\label{pr-among-sa}
  We can for example consider, for every $j\in\{-(r-1)-a,\dotsc,r-1\}$, the projection $\pr_j\colon\Lambda_{\Z_{(p)}}^\S\to \Lambda_{\Z_{(p)}}^j$.

  Since the $G$-action on $\Lambda_{\Z_{(p)}}^\S=\prod_{i=-(r-1)-a}^{r-1}\Lambda_{\Z_{(p)}}^i$ is diagonal, the $\pr_j$ are fixed, i.e., we may count the $\left(\pr_j\right)_j$ among the $\left(s_\alpha\right)_\alpha$.
\end{Remark}

\begin{Lemma}
  $G_{\Z_{(p)}}\otimes_{\Z_{(p)}}\Z_p={\cal G}$.
\end{Lemma}

\begin{Proof}
  First, $G_{\Z_{(p)}}\otimes_{\Z_{(p)}}\Z_p$ is the Zariski closure $G_{\Z_p}$ of $G$ in $\GL(\Lambda^\S_{\Z_p})$ and $G_{\Z_p}\otimes\Q_p=G$: This follows from \cite[Lemma~14.6]{gw}.
  
  According to \stacks{056B}, $G_{\Z_p}$ is the topological closure of $G$ in $\GL(\Lambda^\S_{\Z_p})$ endowed with the reduced subscheme structure. The special fiber therefore is reduced as well, i.e. ($\F_p$ being perfect), geometrically reduced, i.e. (being a group scheme), smooth.
  $G_{\Z_p}$ is flat over $\Z_p$ by \cite[Prop.~14.14]{gw}.
  $G_{\Z_p}$ obviously is of finite presentation over $\Z_p$.
  Hence, $G_{\Z_p}$ is smooth over $\Z_p$. It also is affine and is the stabilizer of $\Lambda^\S_{\Z_p}$ (or rather, $G_{\Z_p}(\Z_p)\subseteq G(\Q_p)$ is (and analogously for $G_{\Z_p}(\bZ_p)\subseteq G(\bQ_p)$)).
\end{Proof}

\begin{Remarks}
  \begin{enumerate}[(1)]
  \item Obtaining an embedding of the Bruhat-Tits group scheme to which we can apply Proposition~\ref{hodge-tensors-kisin} is the main reason for altering the Hodge embedding in the way outlined in Section~\ref{sec:alt-hodge}.
  \item We can also interpret the $s_\alpha$ as tensors in $((\Lambda^\S_{\Z_{(p)}})^*)^\otimes$ and that $G_{\Z_{(p)}}$ also is the Zariski closure of $G$ in $\GL((\Lambda^\S_{\Z_{(p)}})^*)$ (using the contragredient representation).
  \end{enumerate}
\end{Remarks}

\paragraph{Hodge tensors: generic fiber.}
\label{sec:hodge-tens-gener}

\begin{Notation}
  Recall the universal abelian scheme $h\colon \mathcal{A}\to\Sh_K(G,X)$ from Remark~\ref{univ-ab-sch-rat}.

We consider the local system $\mathcal{L}_B:=R^1h_*^\mathrm{an}\underline\Q$ on $\Sh_K(G,X)_\C^\mathrm{an}$ and the flat vector bundle $\mathcal{V}_\mathrm{dR}:=R^1h_*\Omega^\bullet$ with Gau\ss{}-Manin connection $\nabla$.
\end{Notation}

\begin{Lemma} \textnormal{(See \cite[Lemma~2.3.1]{carsch}.)}
  $\mathcal{L}_B$ can be identified with the local system of $\Q$-vector spaces on $\Sh_K(G,X)^\mathrm{an}$ given by the $G(\Q)$-representation $V$ and the $G(\Q)$-torsor
  \begin{equation*}
    p\colon X\times (G(\A_f)/K) \to G(\Q) \backslash (X\times (G(\A_f)/K)=\Sh_K(G,X)(\C).
  \end{equation*}

  Said torsor is isomorphic to the $G(\Q)$-torsor $I:=\underline{\mathrm{Isom}}((V,s_\alpha),({\cal L}_B,s_{\alpha,B}))$ that maps an open subset $U\subseteq \Sh_K(\C)$ to
  \begin{equation*}
    I(U)=\{\beta\colon V\times U\xrightarrow{\sim} \left.{\cal L}_B\right|_U\suchthat \beta(s_\alpha)=s_{\alpha,B}\},
  \end{equation*}
  where the $\left(s_{\alpha,B}\right)_\alpha\subset \mathcal{L}_B^\otimes$ are the global sections corresponding to the $G(\Q)$-invariant tensors $\left(s_\alpha\right)_\alpha$.
\end{Lemma}

\begin{Remark}
  By de Rham's theorem,
  \begin{equation*}
    \mathcal{V}_{\mathrm{dR},\C}^\mathrm{an} \cong \mathcal{L}_B\otimes_\Q\mathcal{O}_{\Sh_K(G,X)^\mathrm{an}}.
  \end{equation*}

  In particular the global sections $\left(s_{\alpha,B}\right)_\alpha\subset \mathcal{L}_B^\otimes$ yield flat global sections $\left(s_{\alpha,\mathrm{dR}}\right)_\alpha\subset (\mathcal{V}_{\mathrm{dR},\C}^\mathrm{an})^\otimes$. All such sections arise from sections in $\mathcal{V}_{\mathrm{dR},\C}^\otimes$ (which we will denote the same), see \cite[\nopp 2.2]{kisin}.
\end{Remark}

\begin{Remark}
  Let $k/E$ be a field extension embeddable into $\C$, and choose an embedding $\Q_p\hookrightarrow\C$ and an $E$-embedding $\sigma\colon\bar k\hookrightarrow\C$. Let $x\in \Sh_K(G,X)(k)$. By $p$-adic Hodge theory, the embedding $\sigma$ gives rise to isomorphisms
  \begin{equation*}
    H^1_\mathrm{dR}(\mathcal{A}_x/k)\otimes_k\C  \xrightarrow{\sim} H^1_B(\mathcal{A}_x(\C),\Q)\otimes_\Q\C \xrightarrow{\sim} H^1_\mathrm{\acute{e}t}(\mathcal{A}_{x,\bar k},\Q_p)\otimes_{\Q_p}\C,
  \end{equation*}
  and by proper base change $(\mathcal{L}_B)_x:=x^{-1}\mathcal{L}_B\cong H^1_B(\mathcal{A}_x(\C),\Q)$.
\end{Remark}

\begin{Notation}
  We denote by $s_{\alpha,B,x}$ the fiber of $s_{\alpha,B}$ at $x$, and by $s_{\alpha,\mathrm{dR},x}$ and $s_{\alpha,\mathrm{\acute{e}t},x}$, respectively, the images of $s_{\alpha,B,x}$ under the above isomorphisms. For the Betti--Ã©tale comparison, we don't need to go all the way to $\C$ and have $s_{\alpha,\mathrm{\acute{e}t},x}\in H^1_\mathrm{\acute{e}t}(\mathcal{A}_{x,\bar k},\Q_p)$.
\end{Notation}

\begin{Remark}
  $(s_{\alpha,\mathrm{dR},x},s_{\alpha,\mathrm{\acute{e}t},x})$ is an absolute Hodge cycle in the sense of \cite{deligne-hc}, for every $\alpha$.
\end{Remark}

\begin{Lemma}\textnormal{\cite[Lemma~2.2.1]{kisin}}
  The $\Gal(\bar k/k)$-action on $H^1_\mathrm{\acute{e}t}(\mathcal{A}_{x,\bar k},\Q_p)$ fixes every $s_{\alpha,\mathrm{\acute{e}t},x}$ and factors through $G(\Q_p)$. Moreover, $s_{\alpha,\mathrm{dR},x}\in H^1_\mathrm{dR}(\mathcal{A}_x/k)$.

  In particular, $(s_{\alpha,\mathrm{dR},x},s_{\alpha,\mathrm{\acute{e}t},x})$ is independent of the choices made above.
\end{Lemma}

\begin{Corollary}\textnormal{\cite[Lemma~2.2.2]{kisin}}
  $s_{\alpha,\mathrm{dR}}$ is defined over $E$ for all $\alpha$, i.e., $s_{\alpha,\mathrm{dR}}\in\mathcal{V}_{\mathrm{dR},E}^\otimes$.
\end{Corollary}

\paragraph{Hodge tensors: special fiber.}
\label{sec:hodge-tens-special}

\begin{Notation}
  Denote by $\mathscr{A}_\mathrm{Siegel}$ the universal abelian scheme over $\sS_J(\GSp(V^\S),S^{\S,\pm})$, cf. Remark~\ref{paragraph-moduli}, and by $\mathscr{A}$ the pullback of $\mathscr{A}_\mathrm{Siegel}$ to $\sS_K(G,X)$. Then the pullback of $\sA$ to $\Sh_K(G,X)_E$ agrees with the pullback of $\mathcal{A}$ of Remark~\ref{univ-ab-sch-rat}.

  We will occasionally call $\sA$ the \defn{``universal'' abelian scheme} (with quotation marks).
\end{Notation}

Let $\bar x\in \mathscr{S}_K(G,X)(\bar\F_p)$ and let $\D_{\bar x}:=\D(\sA_{\bar x}[p^\infty])(W)$ be the DieudonnÃ© module of the associated fiber of the ``universal'' abelian scheme, $W=W(\bar\F_p)=\bZ_p$. Choose a closed point generization $x\in\Sh_K(G,X)(L)$, $L/E$ finite field extension, in the sense of Definition~\ref{def-cpg}. Note that $L$ can be embedded into $\C_p\cong\C$. Then
\begin{equation*}
  H^1_\mathrm{\acute{e}t}(\sA_{x,\bar E},\Z_p)\cong (T_p\sA_{x,\bar E})^*\cong (T_p\sA_{x,\bar E}^*)(-1)\cong (\Lambda^\S)^*,
\end{equation*}
allowing us to identify the tensors $s_\alpha$ with tensors $s_{\alpha,\mathrm{\acute{e}t},x}\in H^1(\sA_{x,\bar E},\Z_p)^\otimes$. Again, these tensors can be shown to be $\Gal(\bar E/L)$-invariant. Now we have the $p$-adic comparison isomorphism
\begin{equation}\label{eq:etcriscomp}
  H_\mathrm{\acute{e}t}^1(A_{x,\bar E},\Q_p)\otimes_{\Q_p}B_\mathrm{cris} \cong H^1_\mathrm{cris}(\sA_{\bar x}/W)\otimes_{W}B_\mathrm{cris}=\D(\sA_{\bar x}[p^\infty])(W)\otimes_{W}B_\mathrm{cris},
\end{equation}
and by \cite[\nopp 3.3.8]{kisin-pappas}, via this isomorphism, the $s_{\alpha,\mathrm{\acute{e}t},x}$ also correspond to tensors $s_{\alpha,0}:=s_{\alpha,0,\bar x}$ in $\D_{\bar x}^\otimes$. In fact, we get an isomorphism
\begin{equation}
  (\Lambda^\S)^*\otimes_{\Z_p}\bZ_p\cong\D_{\bar x}\label{eq:ddmodultriv}
\end{equation}
identifying $s_{\alpha}\otimes 1$ with $s_{\alpha,0,\bar x}$.

\paragraph{Globalizing crystalline tensors.}
\label{sec:glob-cryst-tens}

Above we constructed crystalline tensors fiberwise. This will not suffice to understand the local geometry of the special fiber of the integral model. Therefore we now ``globalize'' the tensors.

Let $\bar x\in\sS_K(G,X)(\bar\F_p)$.
As in Section~\ref{sec:norm-compl}, we write $R_G:=R_{G,\bar x}:=\widehat{\mathcal{O}}_{\sS_K(G,X)_\bOE,\bar x}$.  (So $\Spf R_G$ is the $\bar x$-adic completion of $\sS_K(G,X)_\bOE$.)

\begin{Remark}
  In \cite[\nopp 3.2.12, 3.2.14]{kisin-pappas} a construction of a ``universal'' deformation $p$-divisible group
  \begin{equation}\label{eq:univdeformRG}
    \mathscr{G}_{R_G}\to\Spec R_G
  \end{equation}
  is given, characterized by its DieudonnÃ© display, to wit $\mathrm{DDisp}(\mathscr{G}_{R_G})=\D_{\bar x}\otimes_W\hat{W}(R_G)$ endowed with a natural DieudonnÃ© display structure. 
  
  $\mathscr{G}_{R_G}$ can be identified with the pullback of the $p$-divisible group of the ``universal'' abelian scheme $\sA$.
\end{Remark}

\begin{Remark}
  $R_G$ is a complete local normal noetherian ring with perfect residue field.
\end{Remark}

\begin{Remark}
  By \cite[Theorem~B]{laurel}, the DieudonnÃ© crystal associated with $\mathrm{DDisp}(\sG_{R_G})$ coincides with the DieudonnÃ© crystal $\D(\sG_{R_G})$ of $\sG_{R_G}$.

  In fact, $\mathrm{DDisp}(\sG_{R_G})=\D(\sG_{R_G})^\vee(\hat{W}(R_G))$ endowed with a natural DieudonnÃ© display structure, cf. \cite[\nopp 3.1.7]{kisin-pappas}.

  In particular we get tensors $t_{\alpha,\bar x}^\mathrm{def}\in \D(\sG_{R_G})^\vee(\hat{W}(R_G))^\otimes$ corresponding to $s_{\alpha,0,\bar x}\otimes 1\in (\D_{\bar x}\otimes_W\hat{W}(R_G))^\otimes$, and accordingly tensors $t_{\alpha,\bar x}^\mathrm{def}\otimes 1\in \D(\sG_{R_G})^\vee(W(R_G))^\otimes$.
\end{Remark}

\begin{Definition} \textnormal{(See \cite[Definition~2.8]{hamacher}.)}
  Let $\sG\to S$ be a $p$-divisible group over a formally smooth $\F_p$-scheme $S$. Denote by $\D(\sG)$ its contravariant DieudonnÃ© crystal of $\sG$ \cite[Def.~3.3.6]{bbm}, and by $\1:=\D(\underline{\Q_p/\Z_p}_S)$ the unit object in the tensor category of locally free $\mathcal{O}_{S/\Z_p,\mathrm{CRIS}}$-modules. Note that $\D(\sG)$ comes with a Frobenius morphism $\D(\sG)^{(p)}\to\D(\sG)$, and similar for $\1$.
  
  \begin{enumerate}[(1)]
  \item A \defn{tensor} $t$ of $\D(\sG)$ is a morphism $\1\to \D(\sG)^\otimes$ of locally free $\mathcal{O}_{S/\Z_p,\mathrm{CRIS}}$-modules.\footnote{$(\;)^\otimes$ here (and also in point \eqref{item:crys-tt}) is defined as it is defined for $R$-modules in Prop.~\ref{hodge-tensors-kisin}.}
  \item A tensor $t$ of $\D(\sG)$ is called a \defn{crystalline Tate tensor on $\sG$} if it induces a morphism of $F$-isocrystals $\1\to\D(\sG)[\tfrac1p]^\otimes$.\label{item:crys-tt}
  \end{enumerate}
\end{Definition}

\begin{Remarks}
  \begin{enumerate}[(1)]
  \item To give a tensor as in the preceding definition therefore means the following:
   For every $(U,T,i,\delta)$ with
    \begin{enumerate}[(a)]
    \item an $S$-scheme $U$,
    \item a $\Z_p$-scheme $T$ on which $p$ is locally nilpotent,
    \item a closed $\Z_p$-immersion $i\colon U\hookrightarrow T$,
    \item a pd-structure $\delta$ on the ideal in $\mathcal{O}_T$ defining the immersion $i$, compatible with the canonical pd-structure on $p\Z_p$,
    \end{enumerate}
    functorially to give a morphism $\1(U,T,i,\delta)=\mathcal{O}_T(T)\to \D(\sG)(U,T,i,\delta)^\otimes$ of $\mathcal{O}_T(T)$-modules, i.e., functorially to give elements of $\D(\sG)(U,T,i,\delta)^\otimes$.
  \item \label{item:dJvia} Assume $S=\Spec A$, where $A$ is an $\F_p$-algebra which has a $p$-basis (e.g., $A$ perfect) \emph{or} which satisfies \cite[\nopp (1.3.1.1)]{dJvia}, the latter signifying the existence of an ideal $I\subseteq A$ such that
    \begin{itemize}
    \item $A$ is noetherian and $I$-adically complete,
    \item $A$ is formally smooth as a topological $\F_p$-algebra with the $I$-adic topology,
    \item $A/I$ contains a field with a finite $p$-basis and is finitely generated as an algebra over this field.
    \end{itemize}
    
    Also fix a \emph{lift} $\tilde{A}$ of $A$ in the sense of \cite[Def.~1.2.1]{dJvia}, i.e., a $p$-adically complete $\Z_p$-flat ring $\tilde{A}$ together with an isomorphism $\tilde{A}/p\tilde{A}\cong A$ and a ring endomorphism $\sigma\colon\tilde{A}\to\tilde{A}$ such that $\sigma(a)\equiv a^p\mod{p\tilde{A}}$. Note that if $A$ is perfect, then $\tilde{A}:=W(A)$ with the usual Frobenius lift works.

    Then by \cite[Prop.~1.3.3]{bm-iii} and \cite[Cor.~2.2.3]{dJvia} the category of crystals of quasi-coherent $\mathcal{O}_{\Spec(A)/\Z_p,\mathrm{CRIS}}$-modules is equivalent to the category of $p$-adically complete $\tilde{A}$-modules $M$ endowed with an integrable topologically quasi-nilpotent connection $\nabla\colon M\to M\otimes_{\tilde{A}}\hat{\Omega}_{\tilde{A}}^1$.
  \item In the setting of \eqref{item:dJvia}, say $(M,\nabla)$ is the image of $\D(\sG)$ under the equivalence.  Then $M=\D(\sG)(\tilde{A})$, and to give a tensor $\1\to \D(\sG)^\otimes$ means to give a horizontal section of $M^\otimes$.
  \end{enumerate}
\end{Remarks}

\begin{Proposition} \textnormal{\cite[Prop.~3.3.1, Cor.~3.3.7]{hk}}\label{hk-tensors}
  Denote by $\sG_{\osS^\mathrm{perf}}$ and $\sG_{\hsS}$ the pullbacks of the $p$-divisible group $\sA[p^\infty]$ to the perfection\footnote{This being the inverse perfection of $\sS_K(G,X)\otimes_{\mathcal{O}_E}\bar\F_p$ in the terminology of \cite[Section~5]{perfection}.} $\osS^\mathrm{perf}:=(\sS_K(G,X)\otimes_{\mathcal{O}_E}\bar\F_p)^\mathrm{perf}$ and  the $p$-adic completion $\hsS$ of $\sS_K(G,X)\otimes_{\mathcal{O}_E}\bOE$, respectively.

  For a $p$-divisible group $\sX$ over a $p$-adic formal scheme $\fS$ denote by $P(\sX)$ the locally free $W(\mathcal{O}_\fS)$-module given by $P(\sX)(\Spf A)=\D(\sX_A)^\vee(W(A))$ for every open affine formal subscheme $\Spf A\subseteq \fS$.

  Then there exist tensors $\left(t_\alpha\right)_\alpha\subset P(\sG_{\hsS})^\otimes$ whose pullback to $P(\sG_{R_G,\bar x})^\otimes$ coincides with $t_{\alpha,\bar x}^\mathrm{def}\otimes 1$ for all $\bar x\in\sS_K(G,X)(\bar\F_p)$.
  
  Via pullback, these $t_\alpha$ yield crystalline Tate tensors $(\bar{t}_\alpha)_\alpha$ on $\sG_{\osS^\mathrm{perf}}$.

  Moreover, if $x\in\sS_K(G,X)(\mathcal{O}_L)$, $L/\bE$ finite field extension, is a closed point generization of $\bar x\in\sS_K(G,X)(\bar\F_p)$, then $s_{\alpha,\mathrm{\acute{e}t},x}\in (T_p\sA_{x,\bar{E}})^\otimes$ gets identified with $\bar{t}_{\alpha,\bar x}\in \D(\sA_{\bar x}[p^\infty])(W(\bar\F_p))$ by the $p$-adic comparison isomorphism \eqref{eq:etcriscomp}.
\end{Proposition}

\section{Central leaves in the case of parahoric reduction}
\label{part:central-leaves}

We begin by giving some history. The foliation given by the central leaves was introduced by Oort \cite{oort,oort-dim} in the setting of a $p$-divisible group over a characteristic $p$ scheme (cf. Remark~\ref{rmk-def-leaves} below). As already indicated in the introduction, the idea is to consider for every isomorphism class of $p$-divisible groups over an algebraically closed field of characteristic $p$ the locus where the isomorphism class of the geometric fibers of the $p$-divisible group is the given one. It is not at all readily apparent why this would be a reasonable notion in geometrical or topological terms. The key tool here is given by the slope filtrations introduced by Zink \cite{zink-slope}, and a key insight is that there is a number $N\geq 0$ such that a $p$-divisible group over an algebraically closed field of characteristic $p$ is in the prescribed isomorphism class if and only if the analogous statement holds for the $p^N$-torsion subgroups. For Shimura varieties of Hodge type and good reduction, the central leaves were studied by Mantovan \cite{mantovan-unitary,mantovan} (PEL cases), Vasiu \cite{vasiu-levelm} and Zhang \cite{zhang} (a good survey), among others.

He and Rapoport formulated a set of axioms that  integral models for Shimura varieties with parahoric level are supposed to satisfy in order for them to merit the label ``canonical''. Among these is having a notion of well-behaved central leaves. This is what we investigate in the Hodge type situation. As already mentioned, Zhou \cite{zhou} concurrently and independently also worked on this and obtained very similar results in a different way, and we use some of his results from the first version of the preprint, which didn't yet address the change-of-parahoric map between central leaves. Central leaves in the parahoric level case also are subject in work of Hamacher and Kim \cite{hamacher,kim-leaves,hk}, but they do not deal with changing the parahoric level.

We still fix a Shimura datum $(G,X)$ of Hodge type, a parahoric subgroup $K_p\subseteq G(\Q_p)$ (associated with a Bruhat-Tits group scheme $\mathcal{G}\to\Spec\Z_p$) and a sufficiently small open compact subgroup $K^p\subseteq G(\A_f^p)$. We also keep up our standard assumptions~\ref{std-assum}.

\subsection{Definition of central leaves}
\label{sec:defin-centr-leav}

We use the notation for Hodge tensors, in particular $\left(s_\alpha\right)_\alpha$ and $\left(s_{\alpha,0}\right)_\alpha$, established in Section~\ref{sec:story-integral-model}.

\begin{Definition}
  Two points $\bar{x}_1,\bar{x}_2\in\mathscr{S}_K(G,X)(\bar\F_p)$ lie in the same \defn{central leaf} if there is an isomorphism $\D_{\bar{x}_1}\cong\D_{\bar{x}_2}$ of DieudonnÃ© modules that takes $s_{\alpha,0,\bar{x}_1}$ to $s_{\alpha,0,\bar{x}_2}$.

  They lie in the same \defn{naive central leaf} if there is an arbitrary isomorphism $\D_{\bar{x}_1}\cong\D_{\bar{x}_2}$ of DieudonnÃ© modules (i.e., an isomorphism $\sA_{\bar{x}_1}[p^\infty]\cong\sA_{\bar{x}_1}[p^\infty]$ of $p$-divisible groups).
\end{Definition}

\begin{Remark}\label{rmk-def-leaves}
  More generally, one can make analogous definitions when given any abelian scheme $\mathcal{A}\to S$ over any $\F_p$-scheme $S$ (such that for every point we are given tensors on the DieudonnÃ© module of the corresponding abelian variety). Given two arbitrary points, one compares them after going to a common algebraically closed extension of the residue fields, cf. \cite{oort}.
\end{Remark}

\begin{Notation}
  Recall that we denote the completion of the maximal unramified extension $\Q_p^\mathrm{ur}$ of $\Q_p$ by $\bQ_p$, and its ring of integers accordingly by $\bZ_p$. We set $\bK:=\bK_p:=\mathcal{G}(\bZ_p)$ and define $\bK_\sigma$ to be the graph of the Frobenius $\sigma\colon \bK\to \bK$. So dividing out the action of $\bK_\sigma$ (which is mostly how $\bK_\sigma$ will make an appearance) means dividing out the action of $\bK$ by $\sigma$-conjugation.
\end{Notation}

\begin{Definition}\label{Upsilon}
  The central leaves are the fibers of the map
  \begin{equation*}
    \Upsilon=\Upsilon_K\colon \mathscr{S}_K(G,X)(\bar\F_p)\to G(\bQ_p)/\bK_\sigma
  \end{equation*}
  given as follows: For $\bar x\in\mathscr{S}_K(G,X)(\bar\F_p)$ there is an isomorphism $\beta\colon V_{\Z_p}^*\otimes_{\Z_p}\bZ_p\cong\D_{\bar x}$ (equation~\eqref{eq:ddmodultriv}) sending $s_{\alpha}\otimes 1$ to $s_{\alpha,0,\bar x}$.  We hence can interpret the Frobenius on $\D_{\bar x}$ as an element of $G(\bQ_p)/\bK_\sigma$, where dividing out $\bK_\sigma$ rids us of the ambiguity introduced by the choice of the isomorphism $\beta$.
\end{Definition}

\begin{Notation}
  We write
  \begin{alignat*}{2}
    G(\bQ_p) &\to{}& \centermathcell{C(G):=G(\bQ_p)/\bK_\sigma} &\to B(G):=G(\bQ_p)/G(\bQ_p)_\sigma, \\
    b &\mapsto{}& \centermathcell{\llbracket b\rrbracket} &\mapsto [b].
  \end{alignat*}
\end{Notation}

\subsubsection{An alternative characterization of the central leaves}
\label{sec:an-altern-defin}

Consider two points $\bar{x}_1,\bar{x}_2\in\mathscr{S}_K(G,X)(\bar\F_p)$ in the same central leaf, i.e., such that there is an isomorphism $\D_{\bar{x}_1}\cong\D_{\bar{x}_2}$ of DieudonnÃ© modules that takes $s_{\alpha,0,\bar{x}_1}$ to $s_{\alpha,0,\bar{x}_2}$.

$\bar{x}_j$ for $j=1,2$ has an associated isogeny chain of abelian schemes in $\mathscr{S}_{N}(\GSp(V),S^{\pm})$, cf. Remark~\ref{rz-moduli}, and $\D_{\bar{x}_j}=\D(\prod_{i=-(r-1)-a}^{r-1}A_{j,i})=\prod_{i=-(r-1)-a}^{r-1}\D(A_{j,i})$, where $(A_{j,i})_i$ is the isogeny chain associated with $\bar{x}_j$ under $\mathscr{S}_K(G,X)\to\mathscr{S}_J(\GSp(V^\S),S^{\S,\pm})$, cf. Remark~\ref{paragraph-moduli} and Proposition~\ref{diag-emb}. Note that $A_j=\prod_{i=-(r-1)-a}^{r-1}A_{j,i}$ (plus extra structure) is the point associated with $\bar{x}_j$ under $\mathscr{S}_K(G,X)\to\mathscr{S}_J(\GSp(V^\S),S^{\S,\pm})$.

\begin{Lemma}
  Any tensor-preserving isomorphism $\D_{\bar{x}_1}\cong\D_{\bar{x}_2}$ yields, for all $i$, isomorphisms ${\D(A_{1,i})\cong\D(A_{2,i})}$.
\end{Lemma}

\begin{Proof}
  This follows immediately from Remark~\ref{pr-among-sa}.
\end{Proof}

We may thus rephrase our definition of central leaves as follows.

\begin{Lemma}\label{old-and-new-leaves}
  Let $\bar{x}_1,\bar{x}_2\in\mathscr{S}_K(G,X)(\bar\F_p)$ and, as above, denote by $A_{j,i}$ the abelian varieties in the associated isogeny chains, $j=1,2$.

  $\bar{x}_1$ and $\bar{x}_2$ lie in the same central leaf if and only if there is an isomorphism between the associated rational DieudonnÃ© modules of $A_{j,i}$ (which is common to all $A_{j,i}$, $j$ fixed, $i$ variable) respecting the tensors and identifying the lattices $\D(A_{1,i})$ and $\D(A_{2,i})$ for all $i$.

  ``Respecting the tensors'' means the following: We get an induced identification of the rational DieudonnÃ© modules of $A_1$ and $A_2$. This identification has to preserve the tensors.
\end{Lemma}

\subsection{Local closedness of central leaves}
\label{sec:local-clos-centr}

\subsubsection{Topological lemmas}
\label{sec:topological-lemmas}

We begin with some purely topological preliminaries which we shall use to globalize statements about formal neighborhoods.

\begin{Lemma}
  Let $X$ be a topological space with a subset $A\subseteq X$, and for $x\in X$ let $U_x\subseteq X$ be the set of generizations of $x$. Consider the statements
  \begin{enumerate}[(i)]
  \item $A\cap U_x\subseteq U_x$ closed (resp. open) for all $x\in X$.\label{item:top-i}
  \item $A\cap U_x\subseteq U_x$ closed (resp. open) for all closed points $x\in X$.\label{item:top-ii}
  \item $A\subseteq X$ stable under specialization (resp. generization).\label{item:top-iii}
  \end{enumerate}

  We have \eqref{item:top-i}$\implies$\eqref{item:top-ii},\eqref{item:top-iii}, and if every $x\in X$ has a specialization in $x$ that is a closed point,  we also have \eqref{item:top-ii}$\implies$\eqref{item:top-iii}.
\end{Lemma}

\begin{Proof}
  Denote the closure of $A$ in $X$ by $\cl_X(A)$.

  For ``\eqref{item:top-i}$\implies$\eqref{item:top-iii}'' (with ``closed'' and ``specialization'', respectively) let $x\in A$ and $s\in \cl_X(\{x\})$. Then $s,x\in U_s$, and $s\in\cl_X(\{x\})\cap U_s=\cl_{U_s}(\{x\})$. By assumption $A\cap U_s\subseteq U_s$ is stable under specialization, hence $s\in A\cap U_s\subseteq A$.

  The rest is proven along similar lines.
\end{Proof}

\begin{Example}
  $\N\subseteq\A^1_\C$ satisfies property \eqref{item:top-i} from the lemma but is not closed (hence not constructible).
\end{Example}

\begin{Lemma}\label{top-lemma}
  Let $X$ be a sober topological space and let $A\subseteq X$ be stable under both generization and specialization.
  
  Then $A$ is a union of irreducible components of $X$.

  If additionally $X$ only has finitely many irreducible components, $A$ even is a union of connected components and is open and closed.
\end{Lemma}

\begin{Proof}
  Consider the unique generic points $\left\{\eta_i\right\}_{i\in I}$ of the irreducible components of $X$. By assumption, we have, firstly, that if $A$ contains $\eta_i$ then it contains the entire irreducible component $\overline{\{\eta_i\}}$, and, secondly, that if $A$ contains any one point of $\overline{\{\eta_i\}}$, then it also contains $\eta_i$. Similarly, if $A$ contains any one irreducible component $\overline{\{\eta_i\}}$, it also contains all irreducible components that meet $\overline{\{\eta_i\}}$.
\end{Proof}

\subsubsection{Local closedness}
\label{sec:local-closedness}

To show the local closedness of central leaves, we can content ourselves with a construction in the spirit of Proposition~\ref{hk-tensors} but somewhat simpler. Namely, we pull back \eqref{eq:univdeformRG} to obtain
\begin{align*}
  \mathscr{G}_{\bar x} \to \Spec \bar R_G&:=\Spec R_G\otimes_{\bOE}\bar\F_p\\
                                         &=\Spec (\widehat{\cal O}_{\mathscr{S}_K(G,X)_{\bOE},\bar x}\otimes\bar\F_p)=\Spec(\widehat{\cal O}_{\mathscr{S}_K(G,X)\otimes_{{\cal O}_{\bfE,(v)}}\bar\F_p,\bar x}).
\end{align*}

The DieudonnÃ© display of $\mathscr{G}_{\bar x}$ then is
\begin{equation*}
  \mathrm{DDisp}(\mathscr{G}_{\bar x})=\D(\sG_{\bar x})^\vee(\hat{W}(\bar{R}_G)) =\D_{\bar x}\otimes_{\bZ_p}\hat{W}(\widehat{\cal O}_{\mathscr{S}_K(G,X)\otimes_{{\cal O}_{\bfE,(v)}}\bar\F_p,\bar x}).
\end{equation*}
In particular we can consider the tensors $\left(s_{\alpha,0}\otimes 1\right)_\alpha$ on this DieudonnÃ© display and on $\D(\sG_{\bar x})^\vee(W(\bar{R}_G))$.

Thus we get crystalline Tate tensors $\left(u_{\alpha,\bar x}\right)_\alpha$ on $\D(\sG_{\bar x,\mathrm{perf}})$, where $\sG_{\bar x,\mathrm{perf}}$ is the pullback of $\sG_{\bar x}$ to $\bar{R}_G^\mathrm{perf}$.

\begin{Theorem}
  The central leaves on $\mathscr{S}_K(G,X)\otimes\bar\F_p$ are open and closed in the naive central leaves.
\end{Theorem}

\begin{Proof}
  We consider the $p$-divisible group $\sG_{\bar x,\mathrm{perf}}$ over $\widehat{\cal O}_{\mathscr{S}_K(G,X)\otimes_{{\cal O}_{\bfE,(v)}}\bar\F_p,\bar x}^\mathrm{perf}$ together with its crystalline Tate tensors. Note that perfection makes no difference for topological considerations.

  On the naive central leaves on $\Spec(\widehat{\cal O}_{\mathscr{S}_K(G,X)\otimes_{{\cal O}_{\bfE,(v)}}\bar\F_p,\bar x})=\Spec(\widehat{\cal O}_{\mathscr{S}_K(G,X)\otimes_{{\cal O}_{\bfE,(v)}}\bar\F_p,\bar x}^\mathrm{perf})$ the $p$-divisible group is geometrically fiberwise constant. By a lemma of Hamacher \cite[\nopp 2.12]{hamacher}, the tensors then are geometrically fiberwise constant as well.
  
  Using the fact that $\widehat{\cal O}_{\mathscr{S}_K(G,X)\otimes_{{\cal O}_{\bfE,(v)}}\bar\F_p,\bar x}$ is noetherian, we obtain that the central leaf on $\widehat{\cal O}_{\mathscr{S}_K(G,X)\otimes_{{\cal O}_{\bfE,(v)}}\bar\F_p,\bar x}$ is closed and open in the corresponding naive central leaf.

  $\Spec{\cal O}_{\mathscr{S}_K(G,X)\otimes_{{\cal O}_{\bfE,(v)}}\bar\F_p,\bar x}$ has the quotient topology with respect to the fpqc covering $\Spec\widehat{\cal O}_{\mathscr{S}_K(G,X)\otimes_{{\cal O}_{\bfE,(v)}}\bar\F_p,\bar x}\to \Spec{\cal O}_{\mathscr{S}_K(G,X)\otimes_{{\cal O}_{\bfE,(v)}}\bar\F_p,\bar x}$,
  hence we also get the same result after leaving away completion.

By Lemma~\ref{top-lemma} we conclude that we also get the same result for the central leaves and naive central leaves on $\mathscr{S}_K(G,X)\otimes\bar\F_p$. 
\end{Proof}

\begin{Corollary}\label{leaves-loc-closed}
  The central leaves on $\mathscr{S}_K(G,X)\otimes\bar\F_p$ are locally closed.
\end{Corollary}

\begin{Proof}
  Because Oort \cite{oort} has proven that the naive central leaves on ${\mathscr{S}_K(G,X)\otimes\bar\F_p}$ are locally closed, this is an immediate consequence of the preceding theorem.
\end{Proof}

\begin{Corollary}\label{leaf-in-newton}
  In fact, Oort has shown that the naive central leaf is closed in the naive Newton stratum. So our argument even goes to show that the central leaves are closed in their respective Newton strata (naive and non-naive).
\end{Corollary}

\subsection{Quasi-isogenies of \texorpdfstring{$p$}{p}-divisible groups}
\label{sec:quasi-isogenies-p}

Let $b\in G(\bQ_p)=\mathcal{G}(\bQ_p)\subseteq \GL(\Lambda^\S)(\bQ_p)$ and denote by $b^\S$ the image in $\GL(V^\S)(\bQ_p)$.

Let $J_b$ denote the $\Q_p$-reductive group given on $R$-valued points, $R \in (\Q_p\text{-alg})$, by
\begin{equation*}
  J_b(R) := \left\{ g\in G(R\otimes_{\Q_p}\bQ_p)=\Res_{\bQ_p/\Q_p}(G)(R) \suchthat gb\sigma(g)^{-1} = b \right\}
\end{equation*}
(cf. \cite[Prop.~1.12]{rz} and \cite[Prop.~2.2.6]{kim-leaves}).

$J_b(\Q_p)$ then naturally attains the structure of a locally profinite group. We make it into a formal group scheme $\underline{J_b(\Q_p)}$ over $\Spf\bZ_p$; with $\underline{J_b(\Q_p)}(U)=\mathrm{Map}_\mathrm{cont}(U,J_b(\Q_p))$ for formal test schemes $U\to\Spf\bZ_p$.\footnote{Locally: If $G=\varprojlim_n G_n$ is a profinite set (or even profinite group), then $\underline{G}:=\varprojlim_n \underline{G_n}$ with $\underline{G_n}$ constant formal (group) scheme.}

Let $\X_b=\X_{b,K}=\X_{b^\S}$ be a polarized $p$-divisible group over $\bar\F_p$ with a distinguished isomorphism between its DieudonnÃ© (symplectic) module and $(\Lambda^\S)_{\bZ_p}$ under which the Frobenius is identified with $b$. The existence of such a $p$-divisible group depends on $b$, but we \emph{assume} that $\X_b$ exists from now on (see also \ref{notation-cop}\,\eqref{item:notation-cop-2} below). By DieudonnÃ© theory we obtain a bijection
\begin{align*}
  J_{b^\S}(\Q_p) &:= \left\{ g\in \GL(\Lambda^\S\otimes \bQ_p) \suchthat gb^\S\sigma(g)^{-1} = b^\S \right\} \\
  &\cong \{\text{(self-)quasi-isogenies of } \X_b \text{ (without polarization)}\}
\end{align*}
(and similar with polarization (replace $\GL$ by $\mathrm{(G)Sp}$)), and $J_b(\Q_p)\subseteq J_{b^\S}(\Q_p)$ under this bijection corresponds to the \defn{tensor-preserving} quasi-isogenies of $\X_b$.

We need to understand more than this; at least we will need to understand the tensor-preserving quasi-isogenies of $(\X_b)_{\Omega}$ for all algebraically closed fields of characteristic $p$. Using internal hom $p$-divisible groups first defined by Chai and Oort, Caraiani and Scholze \cite[Prop.~4.2.11]{carsch} worked out the structure of the quasi-isogeny group (albeit in a more special setting than ours; Kim \cite{kim-leaves} generalized it to our setting).

First of all, following \cite[section~3.2]{kim-leaves}, we denote by $\Qisg(\X_b)$ the formal group scheme over $\Spf\bZ_p$ given by
\begin{equation*}
  \Nilp_{\bZ_p} \to \mathrm{Grp}, \quad R \mapsto \{\text{quasi-isogenies of } (\X_b)_{R/p}\},
\end{equation*}
where $\Nilp_{\bZ_p}$ is the category of $\bZ_p$-algebras on which $p$ is nilpotent (anti-equivalent to the category of affine schemes over $\Spf\bZ_p$\footnote{That is, affine schemes together with a natural transformation of functors $\mathrm{Ring}\to \mathrm{Set}$ from it to $\Spf\bZ_p$.}).

\begin{Definition}
  Let $R$ be a topological ring.
  \begin{enumerate}[(1)]
  \item $R$ is called \defn{adic}, if there is an ideal $I$ (called \defn{ideal of definition}) such that $\left\{I^n\right\}_{n\in\N}$ is a basis of neighborhoods (or, equivalently, basis of \emph{open} neighborhoods) of $0$.
  \item $R$ is called \defn{f-adic} (or \defn{Huber ring}), if there exists an open subring $A_0\subseteq A$ that is adic with finitely generated ideal of definition. Such a ring is called a \defn{ring of definition}.
  \end{enumerate}
\end{Definition}

\begin{Definition}
  Let $R$ be a ring of characteristic $p$.
  \begin{enumerate}[(1)]
  \item $R$ is \defn{semiperfect}, if the Frobenius endomorphism $\Phi\colon R\to R$ is surjective.
  \item $R$ is \defn{f-semiperfect}, if it is semiperfect and the perfection $R^\flat:=\varprojlim_{\Phi}R$ (with the inverse limit topology) is f-adic.
  \end{enumerate}
\end{Definition}

There is a functorial construction which assigns to every semiperfect ring $R$ its universal $p$-adically complete pd-thickening $A_\mathrm{cris}(R)$ \cite[Prop.~4.1.3]{SchWein}. Set $B_\mathrm{cris}^+(R):=A_\mathrm{cris}(R)[\frac{1}{p}]$. We have
\begin{equation*}
  \D((\X_{b})_R)(A_\mathrm{cris}(R))\cong \D(\X_b)(\bZ_p)\otimes_{\bZ_p}A_\mathrm{cris}(R) \cong \Lambda^\S\otimes_{\Z_p}A_\mathrm{cris}(R)
\end{equation*}
by construction.

\begin{Assumption}
  We will assume from now on that $[b]\in B(G)$ is \emph{neutral acceptable} in the sense of \cite[Def.~2.3]{rapo-vieh}, i.e., $[b]\in B(G,\{\mu\})$. (Of course this is a completely harmless assumption with regard to the setting of Shimura varieties.)
\end{Assumption}

There exists an internal tensor-preserving quasi-endomorphism $p$-divisible group:

\begin{Lemma} \textnormal{\cite[Lemma~3.1.3]{kim-leaves}}
  There is a $p$-divisible group $\mathcal{H}_b^G$ such that for any f-semiperfect $\bar\F_p$-algebra $R$ there is a natural $\Q_p$-linear isomorphism
  \begin{equation*}
    \tilde{\mathcal{H}}_b^G(R) \cong \End_{(s_\alpha)}((\X_b)_R)\bigl[\tfrac{1}{p}\bigr],
  \end{equation*}
  where $\tilde{\mathcal{H}}_b^G$ is the universal cover of $\mathcal{H}_b^G$, and where on the right hand side we have the set of $\gamma\in\End((\X_b)_R)[\frac{1}{p}]$ such that the endomorphism of $\Lambda^\S\otimes_{\Z_p} B_\mathrm{cris}^+(R)$ induced by $\gamma$ preserves the tensors $s_\alpha\otimes 1$.
\end{Lemma}

What we really need is an internal tensor-preserving quasi-isogeny $p$-divisible group. Still following \cite{kim-leaves}, to this end we consider the group sheaf $\Qisg(\X_b)$ on $\Nilp_{\bZ_p}$ given by $R\mapsto \{\text{self-quasi-isogenies of } (\X_b)_R\}$. This can be realized as a closed formal subscheme of $\tilde{\mathcal{H}}_b^2$ (where $\mathcal{H}_b:=\mathcal{H}_{b^\S}=\mathcal{H}_{b^\S}^{\GL(V^\S)}$ (no tensors)).

The closed formal subscheme $\Qisg_G(\X_b)\subseteq \Qisg(\X_b)$ over $\Spf\bZ_p$ then is defined by
\begin{equation*}
  \Qisg_G(\X_b) = \Qisg(\X_b)\times_{\tilde{\mathcal{H}}_b^2}(\tilde{\mathcal{H}}_b^G)^2.
\end{equation*}

\begin{Proposition} \textnormal{ \cite[Prop~3.2.4]{kim-leaves}}
  We have a natural map $\underline{J_b(\Q_p)}\to\Qisg_G(\X_b)$, which has a natural retraction
  \begin{equation*}
    \Qisg_G(\X_b)\to\underline{J_b(\Q_p)}
  \end{equation*}
  with all fibers isomorphic to $\Spf \bZ_p\llbracket x_1^{p^{-\infty}},\dotsc,x_d^{p^{-\infty}}\rrbracket$ as formal schemes, where $d=\langle 2\rho,\nu_{[b]}\rangle$ with $2\rho$ the sum of all the positive roots of $G_{\bar\Q_p}$ and $\nu_{[b]}$ the dominant\footnote{This requires distinguishing a Weyl chamber as the positive one; this arises from the choice of a convenient maximal torus in \cite[just before Rmk.~2.1.4]{kim-leaves}.} representative of the $G(\bQ_p)$-conjugacy class of $\nu_b\in X_*(G_{\bQ_p})_\Q$.\footnote{So $b\mapsto \nu_{[b]}$ is the \emph{Newton map} in the parlance of \cite[256]{kottwitz-iso2}.}
\end{Proposition}

\begin{Corollary}\label{QisgJb}
  Let $R\in\Nilp_{\bZ_p}$ be such that there is no non-trivial continuous homomorphism of $\bZ_p$-algebras $\bZ_p\llbracket x^{p^{-\infty}}\rrbracket\to R$. Put differently, there is no $r \in R$ having a compatible system of $p$-power roots such that all power series with $\bZ_p$-coefficients  in $r$ and its $p$-power roots converge.

  Then
  \begin{equation*}
    \Qisg_G(\X_b)(R) \cong \underline{J_b(\Q_p)}(R).
  \end{equation*}
\end{Corollary}

\begin{Proof}
  With $d$ as in the proposition, the assumption implies that there is no non-trivial continuous homomorphism of $\bZ_p$-algebras $\bZ_p\llbracket x_1^{p^{-\infty}},\dotsc,x_d^{p^{-\infty}}\rrbracket\to R$.
  
  Define $\Qisg_G^\circ(\X_b):=\ker(\Qisg_G(\X_b)\to\underline{J_b(\Q_p)})\cong \Spf \bZ_p\llbracket x_1^{p^{-\infty}},\dotsc,x_d^{p^{-\infty}}\rrbracket$ and consider the split exact sequence
  \begin{equation*}
    0\to \Qisg_G^\circ(\X_b) \to \Qisg_G(\X_b)\to\underline{J_b(\Q_p)}\to 0.
  \end{equation*}
\end{Proof}

\begin{Example}
  Say $R\in\Nilp_{\bZ_p}$ has characteristic $p$, i.e., $R$ is a $\bar\F_p$-algebra and the $p$-adic topology is the discrete topology. Let $r$ be an element of $R$ having a compatible system of $p$-power roots such that all power series with $\bZ_p$-coefficients in $r$ and its $p$-power roots converge. Then $r$ must be nilpotent. Therefore, if $R\in\Nilp_{\bZ_p}$ has characteristic $p$ and is reduced (e.g., is a field), then the conditions of Corollary~\ref{QisgJb} are satisfied.
\end{Example}

In particular we deduce:
\begin{Example}\label{QisgJbField}
  Let $L/\bar\F_p$ be a field extension. Then
  \begin{equation*}
    \Qisg_G(\X_b)(L) \cong J_b(\Q_p).
  \end{equation*}
\end{Example}

\subsection{Almost product structure}
\label{sec:almost-product}

\begin{Notation}
  We shorten notation by writing
  \begin{equation*}
    \sS:=\sS_K(G,X)\otimes_{\mathcal{O}_{\bfE,(v)}}\bOE
    \quad\text{and}\quad
    \osS:=\sS_K(G,X)\otimes_{\mathcal{O}_{\bfE,(v)}}\bar\F_p,
  \end{equation*}
  as well as
  \begin{equation*}
    \sS^\S:=\sS_J(\GSp(V^\S),S^{\pm,\S})\otimes_{\Z_{(p)}}\bZ_p
    \quad\text{and}\quad
    \osS^\S:={\sS_J(\GSp(V^\S),S^{\pm,\S})\otimes_{\Z_{(p)}}\bar\F_p}.
  \end{equation*}

  Moreover, we denote the central leaf associated with $b$ by $\oC^b:=\oC^b_K:=\Upsilon_K^{-1}(\llbracket b\rrbracket)\subseteq\osS$, and by $\oC^{b^\S}\subseteq\osS^\S$ the corresponding central leaf in $\osS^\S$.
\end{Notation}

\begin{Assumption}\label{assumpt-axiom-a}
  From now on, we assume that Axiom~A from \cite{hk} holds. 
\end{Assumption}

\begin{Remark}\label{a-erkl}
  To explain Axiom~A: Consider the \defn{affine Deligne-Lusztig variety}
  \begin{equation*}
    X_\mu(b)_{K}(\bar\F_p):=\left\{ g\bK_p \in G(\bQ_p)/\bK_p\suchthat g^{-1}b\sigma(g)\in \bK_p \Adm(\{\mu\}) \bK_p\right\}
  \end{equation*}
  (here $\Adm(\{\mu\})$ is the  \defn{admissible subset} as defined in \cite{prs}).

  Choose a point $\bar x\in \osS(\bar\F_p)$ with a quasi-isogeny $j\colon \X_b\to \sA_{\bar x}[p^\infty]$ compatible with the extra structure.
  
  We have a map $i^\S_{(\bar x,j)}\colon X_{\mu^\S}(b^\S)(\bar\F_p)\to\osS^\S(\bar\F_p)$ given as follows: For $g \breve{J}_p\in X_{\mu^\S}(b^\S)(\bar\F_p)$, there is a $p$-divisible group $g\X_b$ isogenous to $\X_b$ with Frobenius on the DieudonnÃ© module given by $g^{-1}b\sigma(g)$. The quasi-isogeny $g\X_b\to \X_b\to \sA_{\bar x}[p^\infty]$ lifts to a quasi-isogeny of abelian schemes\footnote{Idea: A lift of a quasi-isogeny $y\colon A[p^\infty]\to X$, $A$ abelian scheme, $X$ $p$-divisible group, is given by $A \to A/\ker(y)$.}. We also get a polarization and a level structure away from $p$ on the new abelian scheme, which is the image of $g\breve{J}_b$ under $i^\S_{(\bar x,j)}$.

  The content of Axiom~A (= Assumption~\ref{assumpt-axiom-a}) now is that
  \begin{equation*}
    X_{\sigma(\mu)}(b)_K(\bar\F_p) \hookrightarrow X_{\mu^\S}(b^\S)(\bar\F_p)\xrightarrow{i^\S_{(\bar x,j)}} \osS^\S(\bar\F_p)
  \end{equation*}
  is to factor through $\osS^-(\bar\F_p)$ such that there is a unique lift $i\colon X_{\sigma(\mu)}(b)_K(\bar\F_p) \to \osS(\bar\F_p)$ with $s_{\alpha,0,i(g\bK_p)}=s_{\alpha,0,\bar x}$ for all $g\bK_p\in X_{\sigma(\mu)}(b)_K(\bar\F_p)$.
\end{Remark}

\begin{Remark}
  Assumption~\ref{assumpt-axiom-a} holds in the hyperspecial case \cite[Thm.~1.4.4]{kisin-modp}, and it holds if $G$ is residually split \cite[Prop.~6.4]{zhou}.
\end{Remark}

\begin{Remark}
  With Assumption~\ref{assumpt-axiom-a} in place, the central leaf $\oC^b$ of $b$ is non-empty by \cite[Rmk.~4.3.1]{hk}.
\end{Remark}

We review some constructions from the paper \cite{hk}. By $\oIg^{b^\S}$ we denote the (special fiber of the) Igusa variety over $\oC^{b^\S}$, cf. \cite[Prop.~1.12]{carsch}. It parameterizes isomorphisms between the standard $p$-divisible group and the universal one, and it is a perfect scheme over $\bar\F_p$.
Define $\oIg^{b,\diamond}:=\oIg^{b^\S}\times_{\osS^{\S,\mathrm{perf}}}\osS^{\mathrm{perf}}$. Furthermore, $\fIg^{b^\S},\fIg^{b,\diamond}$ are the unique flat lifts\footnote{Locally: Let $R$ be a perfect $\bar\F_p$-algebra. Then $W_{\mathcal{O}_E}(R)$ (ramified Witt vectors, cf. Remark~\ref{ram-witt}) is the unique $\varpi$-adically complete $\bOE$-flat algebra lifting it (cf. \cite[Prop.~1.3.3]{ahsendorf}).}
of $\oIg^{b^\S}$ and $\oIg^{b,\diamond}$ over $\Spf\bOE$, given by the formula
\begin{equation}\label{eq:fIg-constr}
  \fIg(R) := \oIg(R/\varpi) \quad \text{for all } R\in\Nilp_{\bOE}
\end{equation}
in both cases, where $\varpi$ is a uniformizer of $\bOE$.

By definition, we have an isomorphism $j\colon \X_b\times \oIg^{b,\diamond}\xrightarrow{\cong}\sA_{\oC^b}[p^\infty]\times_{\oC^b}\oIg^{b,\diamond}$ and $\oIg^b$ is by definition the locus of geometric points of $\oIg^{b,\diamond}$, where $j$ respects the crystalline Tate tensors. This is a closed union of connected components of $\oIg^{b,\diamond}$ \cite[Def./Lemma~5.1.1]{hk}.

$\Qisg(\X_b)$ acts on $\oIg^{b^\S}$ and $\fIg^{b^\S}$, $\oIg^b\to \oIg^{b^\S}$ is a closed immersion with $\Qisg_G(\X_b)_{\bar\F_p}$-stable image, and $\fIg^b\to\fIg^{b^\S}$ is a closed immersion with $\Qisg_G(\X_b)$-stable image \cite[Prop.~5.1.2, Cor~5.1.3]{hk}.

Let $\fM^{b^\S}\to\Spf\bZ_p$ be the Rapoport-Zink space given by
\begin{equation*}
  \fM^{b^\S}(R) =
  \left\lbrace (\sX,\rho) \;\middle|\;
    \begin{tabular}{@{}l@{}}
      $\sX/R$  polarized  $p$-divisible group and \\
      $\rho\colon \X_b\otimes R/p \to \sX\otimes R/p$ quasi-isogeny resp. polarizations
   \end{tabular}
 \right\rbrace.
\end{equation*}

Choose a point $\bar x\in \osS(\bar\F_p)$ with a quasi-isogeny $j\colon \X_b\to \sA_{\bar x}[p^\infty]$ compatible with the extra structure.

$\fM^{b^\S}$ comes with the Rapoport-Zink uniformization map $\Theta^\S\colon \fM^{b^\S}\to\sS^{\S}$ depending on this choice \cite[\nopp (6.3)]{rz}: Given $(\sX,\rho)$ as above, $\rho$ lifts to a quasi-isogeny $\tilde{\rho}\colon \X_b\otimes R\to \sX$, and there is an abelian scheme $\sY$ with $p$-divisible group $\sX$ and a unique lift $\tilde{\sA}_{\bar x}\otimes R\to \sY$ of $\tilde{\rho}$ for a chosen lift $\tilde{\sA}_{\bar x}$ of $\sA_{\bar x}$ to $\bZ_p$. We also get a polarization and a level structure away from $p$ on $\sY$, and $\sY$ with these extra structures is the image of $(\sX,\rho)$ under $\Theta^\S$.

\begin{Remark}
  The affine Deligne-Lusztig variety is the perfection of the special fiber of the Rapoport-Zink space \cite[Prop.~0.4]{zhuaff}. Under this isomorphism, $\Theta^\S$ corresponds to $i^\S$ from Remark~\ref{a-erkl}.
\end{Remark}

Define $\fM^{b,\diamond}:=\fM^{b^\S}\times_{\sS^\S}\sS$ and define $\Theta^\diamond\colon\fM^{b,\diamond}\to\sS$ to be the base change of $\Theta^\S$.

Let $\sX^\diamond\to \fM^{b,\diamond}$ be the pullback of the universal $p$-divisible group over $\fM^{b^\S}$. Then we have two families of tensors on $\D(\sX^\diamond)[\tfrac1p]$ (or more precisely, two families of maps $\fM^{b,\diamond}(\bar\F_p)\ni \bar y\mapsto$ tensor on $\D(\sX^\diamond_{\bar y})[\tfrac1p]$):
\begin{enumerate}[(1)]
\item $\left(t_\alpha^\diamond\right)_\alpha$ obtained from $\left(s_\alpha\right)_\alpha$ via $\D(\X_b\otimes\fM^{b,\diamond})[\tfrac1p]\cong\D(\sX^\diamond)[\tfrac1p]$, and
\item for every $\bar y\in \fM^{b,\diamond}(\bar\F_p)$, $\left(u_{\alpha,\bar y}^\diamond\right)_\alpha$ with $u_{\alpha,\bar y}^\diamond:=(\Theta^\diamond)^*s_{\alpha,0,\Theta^\diamond(\bar y)}$.
\end{enumerate}

$\fM^b$ is defined to be the formal subscheme of $\fM^{b,\diamond}$ corresponding to the locus where these agree. Also, $\Theta\colon \fM^b\to \sS$ is defined to be the restriction of $\Theta^\diamond$. Moreover, $\fM^b$ has a natural $\Qisg_G(\X_b)$-action.

By \cite[Â§\,4]{carsch}, there is an isomorphism
\begin{equation}\label{eq:siegel-aps}
  \fIg^{b^\S}\times\fM^{b^\S} \xrightarrow{\sim} \fX^{b^\S}
\end{equation}
with the \defn{Newton-Igusa variety}
\begin{equation*}  
  \fX^{b^\S}(R) =
\left\lbrace (A,\lambda,\eta^p;\psi)\;\middle|\;
  \begin{tabular}{@{}l@{}}
    $(A,\lambda,\eta^p)\in \sS^\S(R)$, \\
    $\psi\colon (\X_b,\lambda_{\X_b})\otimes R/p \to (A[p^\infty],\lambda)\otimes R/p$ quasi-isogeny
   \end{tabular}
  \right\rbrace.
\end{equation*}

\begin{Remark}\label{carsch-concrete}
  Let us quickly recall how the isomorphism~\eqref{eq:siegel-aps} works (suppressing polarizations in the notation for simplicity). Let $(\mathcal{A},\xi)\in \fIg^{b^\S}(R)$ and $(\sX,\rho)\in\fM^{b^\S}(R)$ be given. Consider the composition
  \begin{equation}\label{eq:the-compos}
    \sX \xrightarrow{\rho^{-1}} \X_b\otimes R \xrightarrow{\xi} \mathcal{A}[p^\infty],
  \end{equation}
  where we denote a lift $\X_b\otimes R/p \xrightarrow{\rho} \sX\otimes R/p$ by $\rho$ again. Lift \eqref{eq:the-compos} to a quasi-isogeny of abelian schemes $\mathcal{A}'\to \mathcal{A}$. Then $(\mathcal{A}',\rho)\in \fX^{b^\S}(R)$ is our image point.
\end{Remark}

Define $\fX^b\subseteq\fX^{b^\S}$ to be the image of $\fIg^b\times \fM^b$ under the isomorphism \eqref{eq:siegel-aps}. This comes with a natural $\Qisg_G(\X_b)$-action; cf. \cite[section~5.2]{hk}.

We have canonical maps $\pi_\infty^\S\colon \fX^{b^\S}\to \sS^\S$ and $\pi_\infty\colon \fX^b \to \sS$.

\begin{Remark}\label{newton-igusa-moduli}
  By \cite[Thm.~5.2.6~(1)]{hk}, $\pi_{\infty,\bar\F_p}^\mathrm{perf}\colon \fX^{b,\mathrm{perf}}_{\bar\F_p}\to\osS^\mathrm{perf}$ represents the moduli problem
  \begin{align*}
    (\text{perfect affine } \osS^\mathrm{perf}\text{-schemes}) &\to \mathrm{Set}, \\
    (\Spec(R)\xrightarrow{Q}\osS^\mathrm{perf}) &\mapsto
                                                  \begin{Bmatrix}
                                                       \psi\colon (\X_b)_R\to(\sG_{\osS^\mathrm{perf}})_Q \text{ quasi-isogeny} \\ \text{compatible with crystalline Tate-tensors}
                                                     \end{Bmatrix}
  \end{align*}
  with $\sG_{\osS^\mathrm{perf}}$ and crystalline Tate tensors on the right hand side as in Proposition~\ref{hk-tensors}.
\end{Remark}

\subsection{Change of parahoric level}
\label{sec:change-parah-level}

Now we consider the question of how the central leaves behave when the level at $p$ is changed from $K_p$ to a larger $K_p'$, where $K_p$ and $K_p'$ are associated with points $x,x'\in\mathcal{B}(G,\Q_p)$ as described in Chapter~\ref{cha:gener-prel}. Define $K:=K_pK^p$, as usual, and $K':=K_p'K^p$. Let $b\in G(\bQ_p)$.

Let $(\mathcal{L},c)$ and $(\mathcal{L}',c')$ be the graded lattice chains that are the images of $x$ and $x'$, respectively, under the embedding of buildings $\iota\colon\mathcal{B}(G,\Q_p)\hookrightarrow\mathcal{B}(\GSp(V),\Q_p)$ described in Section~\ref{sec:bruh-tits-build}. Note that the gradings $c$ and $c'$ play no role for our purposes.

\begin{Remarks}\label{thinout}
  \begin{enumerate}[(1)]
  \item \label{item:thinout} Observe that $K_p$ and $K_p'$ depend only on the minimal facet $\ff$ and $\ff'$ with $x\in\ff$ and $x'\in\ff'$, respectively. The inclusion $K_p\subseteq K_p'$ means precisely that $\ff'\subseteq\overline{\ff}$. We can move $x'$ arbitrarily close to $x$ without altering $K_p'$. The minimal facets $\fg$ and $\fg'$ of $\mathcal{B}(\GSp(V),\Q_p)$ containing $\iota(x)$ and $\iota(x')$, respectively, do \emph{not} depend only on $\ff,\ff'$. Still, for $x'$ sufficiently close to $x$, we may assume that $\fg'\subseteq\overline{\fg}$ (if it were not so, we would not be able to move $\iota(x')$ arbitrarily close to $\iota(x)$ --- but recall that $\iota$ is continuous), i.e., that $\mathcal{L}'$ is a thinned out version of $\mathcal{L}$. Say $\mathcal{L}$ is given by
    \begin{equation*}
      \Lambda^0 \supsetneqq \Lambda^1 \supsetneqq \dotsb\supsetneqq \Lambda^{r-1} \supsetneqq p\Lambda^0
    \end{equation*}
    (as in \eqref{eq:numbered-lattice-chain}); then $\mathcal{L}'$ is given by $(\Lambda^{i_j})_j$ for some $0\leq i_1 < i_2 < \dotsb < i_s < r$.
  \item  In proofs, we can often reduce to the case where ``$\mathcal{L}$ and $\mathcal{L}'$ differ by one element'', which of course is supposed to mean that $s=r-2$ in this notation.
  \end{enumerate}
\end{Remarks}

\begin{Notation}
We define $N_p,N_p',J_p,J_p',\Lambda^\S,\left.\Lambda'\right.^\S,V^\S,\left.V'\right.^\S$ as in Section~\ref{sec:siegel-integr-model}.
\end{Notation}

\subsubsection{The change-of-parahoric morphism}
\label{sec:change-parah-morph}

As explained in Section~\ref{sec:maps-between-shimura}, we obtain (for $K^p$ sufficiently small and some compact open subgroup $N^p\subseteq\GSp(V)(\A_f^p)$) a commutative diagram
\begin{equation*}
  \begin{tikzpicture}[node distance=4cm, auto]
    \node (X) at (0,3) {$\sS_{K_pK^p}(G,X)$};
    \node (Z) at (8,3) {$\sS_{K'_p\left.K\right.^p}(G,X)$};
    \node (Xs) at (0,0) {$\sS_{N_pN^p}(\GSp(V),S^\pm)_{\mathcal{O}_{\bfE,(v)}}$};
    \node (Zs) at (8,0) {$\sS_{N_p'\left.N\right.^p}(\GSp(V),S^\pm)_{\mathcal{O}_{\bfE,(v)}}$};
    
    \draw[->] (X) to node[swap] {finite} (Xs);
    \draw[->] (Z) to node {finite} (Zs);
    \draw[->] (X) to node {$\pi_{K_pK^p,K_p'\left.K\right.^p}$} (Z);
    \draw[->] (Xs) to node {$\pi_{N_pN^p,N_p'\left.N\right.^p}$} (Zs);
  \end{tikzpicture}
\end{equation*}

\begin{Notation}\label{notation-cop}
  \begin{enumerate}[(1)]
  \item The change-of-parahoric map $\pi_{K,K'}:=\pi_{K_pK^p,K_p'\left.K\right.^p}$ restricts to a change-of-parahoric map $\Upsilon_{K_pK^p}^{-1}(\llbracket b\rrbracket)\to\Upsilon_{K_p'K^p}^{-1}(\llbracket b\rrbracket)$ between leaves, which we denote by $\left.\pi_{K,K'}\right|_{\Upsilon_{K}^{-1}(\llbracket b\rrbracket)}$ or simply by $\pi_{K,K'}$ again.
  \item We choose a base point $\bar{x}_b\in\oC^b_K(\bar\F_p)$ and take its image under $\pi_{K,K'}$ as a base point $\bar{x}_b'\in\oC^b_{K'}(\bar\F_p)$. By $\X_{b,K}$ and $\X_{b,K'}$, respectively, we denote the corresponding polarized $p$-divisible groups.\label{item:notation-cop-2}
  \item Denote by $\sA_K$ the universal polarized abelian scheme over $\osS_K$. By slight abuse of notation, we will also use the same notation for its pullback to $\osS_K^\mathrm{perf}$.
  \end{enumerate}
\end{Notation}

\begin{Remark}
  We have morphisms $\X_{b,K}\to \X_{b,K'}$ and, more generally, $\sA_K\to \sA_{K'}$ (lying over $\osS_K\to\osS_{K'}$).  This follows from the ``thinning out'' interpretation, Remark~\ref{thinout}~\eqref{item:thinout}.
\end{Remark}

\subsubsection{Newton-Igusa variety and change-of-parahoric}
\label{sec:newton-igusa-variety}

\begin{Lemma}
  There are natural change-of-parahoric morphisms
  \begin{equation*}
    \fIg^b_K \to \fIg^b_{K'}, \quad \fM^b_K \to \fM^b_{K'}, \quad \fX^b_K \to \fX^b_{K'}.
  \end{equation*}
  lying over the change-of-parahoric morphism $\sS_K\to \sS_{K'}$ and compatible with the Siegel embeddings.
\end{Lemma}

\begin{Proof}
  Consider the universal isomorphism $j\colon \X_{b,K}\times \oIg_K^{b,\diamond}\xrightarrow{\cong} \sA_{K,\oC^b_K}[p^\infty]\times_{\oC^b_K}\oIg_K^{b,\diamond} = \sA_{K}[p^\infty]\times_{\osS_K^\mathrm{perf}}\oIg_K^{b,\diamond}$ and the diagram
  \begin{equation*}
    \begin{tikzpicture}[node distance=4cm, auto,baseline=(current  bounding  box.center)]
      \node (LO) at (0,3) {$\X_{b,K}\times\oIg^{b,\diamond}_K$};
      \node (RO) at (8,3) {$\sA_K[p^\infty]\times_{\osS_K^\mathrm{perf}}\oIg_K^{b,\diamond}$};
      \node (LU) at (0,0) {$\X_{b,K'}\times\oIg^{b,\diamond}_K$};
      \node (RU) at (8,0) {$\sA_{K'}[p^\infty]\times_{\osS_{K'}^\mathrm{perf}}\oIg_K^{b,\diamond}$};
      
      \draw[->] (LO) to  (LU);
      \draw[->] (RO) to  (RU);
      \draw[->] (LO) to node {$\cong$} (RO);
      \draw[dashed,->] (LU) to node{$\cong$} (RU);
    \end{tikzpicture}
  \end{equation*}

  The dashed arrow exists on $\oIg^b_K$, i.e., $\X_{b,K'}\times\oIg^{b}_K\xrightarrow{\cong}\sA_{K'}[p^\infty]\times_{\osS_{K'}^\mathrm{perf}}\oIg_K^{b}$ exists. The reason is that the tensors are respected on $\oIg_K^b$ which in particular forces  $\X_{b,K}\times\oIg^{b,\diamond}_K\xrightarrow{\cong}\sA_K[p^\infty]\times_{\osS_K^\mathrm{perf}}\oIg_K^{b,\diamond}$ to be ``diagonal'', similar to the proof of Proposition~\ref{diag-emb}.

  This isomorphism $\X_{b,K'}\times\oIg^{b}_K\xrightarrow{\cong}\sA_{K'}[p^\infty]\times_{\osS_{K'}^\mathrm{perf}}\oIg_K^{b}$ now corresponds to a morphism $\oIg_K^b\to \oIg_{K'}^b$. By construction \eqref{eq:fIg-constr} of $\fIg^b$, we also get a corresponding morphism $\fIg_K^b\to \fIg_{K'}^b$.

  For $\fX$ and $\fM$, it is very similar.
\end{Proof}

\begin{Lemma}\label{cop-X}
  The isomorphism $\fIg^b\times\fM^b\cong\fX^b$ is compatible with change-of-parahoric, where the change-of-parahoric map for $\fIg^b\times\fM^b$ is by definition the product of those for $\fIg^b$ and $\fM^b$, respectively. 
\end{Lemma}

\begin{Proof}
  All of $\fIg^b,\fM^b,\fX^b$ are embedded into the corresponding Siegel versions as closed formal subschemes, compatible with change-of-parahoric. Hence we may assume that we are in the Siegel case.

  In that case, the isomorphism has a very concrete description, cf. Remark~\ref{carsch-concrete}, from which the lemma is clear.
\end{Proof}

\begin{Definition}
  By $\osS_K^b:=\osS_K^{[b]}$ we denote the Newton stratum associated with ${[b]\in B(G)}$, i.e., $\bar{S}_{K,[b]}$ in the notation of \cite{he-rapo}.
\end{Definition}

\begin{Lemma}
  Let $\Omega/\F_p$ be an algebraically closed field.

  The action of $J_b(\Q_p)$ on the fibers of $\fX_K^b(\Omega)\to\osS^b_K(\Omega)$ is simply transitive and the change-of-parahoric map
  \begin{equation*}
    \fX_K^b(\Omega) \to \fX_{K'}^b(\Omega)
  \end{equation*}
  is $J_b(\Q_p)$-equivariant.
\end{Lemma}

\begin{Proof}
  The moduli description of Remark~\ref{newton-igusa-moduli} gives us that the action is simply transitive upon noting that
  \begin{equation*}
    \Qisg_G(\X_{b,K})(\Omega)\overset{\ref{QisgJbField}}\cong J_b(\Q_p)\overset{\ref{QisgJbField}}\cong\Qisg_G(\X_{b,K'})(\Omega).
  \end{equation*}

  Equivariance follows from the description of the action of $J_b(\Q_p)$ and of the change-of-parahoric map in terms of lattice chains in DieudonnÃ© theory: Given a point $Q\colon\Spec \Omega\to\osS_K^b$, we consider the map between the fiber of $Q$ under $\fX^b_K(\Omega)\to \osS^b_K(\Omega)$ and the fiber of the image of $Q$ under $\fX^b_{K'}(\Omega)\to \osS^b_{K'}(\Omega)$.
  We identify the lattice chain associated with $Q$ with a standard lattice chain such that the Frobenius is identified with a $b'\in G(\bQ_p)$, $b'\equiv b\mod \bK_{p,\sigma}$. Elements of $J_b(\Q_p)$ then act naturally on the common rational DieudonnÃ© module. The action on the fiber of $Q$ is given by altering the quasi-isogenies appearing in the description of that set by the quasi-isogeny obtained this way. Passing from $K$ to $K'$, i.e., from $Q$ to the image of $Q$, means leaving out parts of the lattice chain and enlarging $\bK_p$ (so one still has the same $b'$).
\end{Proof}

\begin{Proposition}\label{cop-X-surj}
  The change-of-parahoric map
  \begin{equation*}
    \fX^b_{K,\bar\F_p} \to \fX^b_{K',\bar\F_p}
  \end{equation*}
  is surjective.
\end{Proposition}

\begin{Proof}
  We check that
  \begin{equation*}
    \fX^b_K(\Omega) \to \fX^b_{K'}(\Omega)
  \end{equation*}
  is surjective for every algebraically closed field $\Omega/\F_p$.

  Considering the diagram
  \begin{equation*}
    \begin{tikzpicture}[node distance=4cm, auto,baseline=(current  bounding  box.center)]
      \node (LO) at (0,3) {$\fX^{b}_{K}(\Omega)$};
      \node (RO) at (3,3) {$\fX^{b}_{K'}(\Omega)$};
      \node (LU) at (0,0) {$\osS_{K}^{b}(\Omega)$};
      \node (RU) at (3,0) {$\osS_{K'}^{b}(\Omega)$};
      
      \draw[->] (LO) to  (LU);
      \draw[->] (RO) to  (RU);
      \draw[->] (LO) to (RO);
      \draw[->>] (LU) to  (RU);
    \end{tikzpicture}
  \end{equation*}
  this follows from the preceding lemma and the fact\footnote{\label{fn:restr-newton}This follows simply from the surjectiveness of $\osS_{K}\to \osS_{K'}$ and the commutativity of the diagram \begin{tikzpicture}[node distance=4cm, auto,baseline=(current  bounding  box.center)]
      \node (LO) at (0,0.7) {$\osS_{K}$};
      \node (RO) at (2,0.7) {$\osS_{K'}$};
      \node (LU) at (1.0,0) {$B(G)$};
      
      \draw[->] (LO) to  (LU);
      \draw[->] (RO) to  (LU);
      \draw[->] (LO) to (RO);
    \end{tikzpicture}, non-horizontal maps being the Newton maps.} that the change-of-parahoric map between Newton strata is surjective.
\end{Proof}

\begin{Corollary}\label{cop-Ig-isom}
  The map $\oIg_K^{b}\to\oIg_{K'}^{b}$ is an isomorphism.
\end{Corollary}

\begin{Proof}
  Lemma~\ref{cop-X} and Proposition~\ref{cop-X-surj} imply that it is surjective.
  
  Also, it is the restriction to closed subschemes of the isomorphism $\oIg_{K}^{b^\S}\to\oIg_{K'}^{b^\S}$. Since all involved schemes are reduced, the corollary follows.
\end{Proof}

\begin{Corollary}\label{cop-surj}
  The change-of-parahoric morphism between central leaves is surjective.
\end{Corollary}

\begin{Proof}
  We have a diagram
  \begin{equation*}
    \begin{tikzpicture}[node distance=4cm, auto]
      \node (LO) at (0,3) {$\oIg_K^{b}$};
      \node (RO) at (3,3) {$\oIg_{K'}^{b}$};
      \node (LU) at (0,0) {$\Upsilon^{-1}_K(b)$};
      \node (RU) at (3,0) {$\Upsilon^{-1}_{K'}(b)$};
      
      \draw[->>] (LO) to  (LU);
      \draw[->>] (RO) to  (RU);
      \draw[->] (LO) to node { $\cong$ }  (RO);
      \draw[->] (LU) to (RU);
    \end{tikzpicture}
  \end{equation*}
  where we already know that all maps but the lower horizontal one are surjective.
\end{Proof}

\begin{Corollary}
  The separable rank of $\Upsilon^{-1}_K(b)\to \Upsilon^{-1}_{K'}(b)$, i.e., the number of geometric points of the  fibers, is finite and constant.
\end{Corollary}

\begin{Proof}
  We have
  \begin{equation*}
    \Upsilon^{-1}_K(b) \cong \oIg_{K}^{b}/\Aut(\X_{b,K}) \cong \oIg_{K'}^{b}/\Aut(\X_{b,K}) \quad\text{and}\quad \Upsilon^{-1}_{K'}(b)\cong \oIg_{K'}^{b}/\Aut(\X_{b,K'}),
\end{equation*}
so that all fibers are isomorphic to $\Aut(\X_{b,K'})/\Aut(\X_{b,K})$.
\end{Proof}

\begin{Corollary}\label{finite-cop}
  The change-of-parahoric map between central leaves is finite.
\end{Corollary}

\begin{Proof}
  We have just seen it to be quasi-finite.

  Corollary~\ref{leaf-in-newton}, combined with the properness of the change-of-parahoric itself and therefore of the change-of-parahoric map restricted to Newton strata (similar argument as in footnote~\ref{fn:restr-newton}), implies the properness of the change-of-parahoric map between central leaves.
\end{Proof}

\begin{Remark}
  By \cite[Cor.~5.3.1]{kim-leaves}, leaves are equidimensional smooth and the dimension of a leaf depends only on the Newton stratum it is in; in particular, if we consider a change-of-parahoric map between leaves,
  \begin{equation*}
    \pi_{K,K'}\colon\Upsilon^{-1}_K(y)\to \Upsilon_{K'}^{-1}(y'),
  \end{equation*}
  then $\dim\Upsilon^{-1}_K(y)=\dim \Upsilon_{K'}^{-1}(y')$.
\end{Remark}

\begin{Corollary}
  The change-of-parahoric map between central leaves  is finite locally free.
\end{Corollary}

\begin{Proof}
  Combine Corollary~\ref{finite-cop} with the preceding remark and \cite[Cor.~14.128]{gw}\footnote{Note that ``$y\in Y$'' may be replaced by ``$y\in f(X)$'' in the statement of \cite[Cor.~14.128]{gw}.}.
\end{Proof}

\begin{Corollary}\label{cop-comp}
  The change-of-parahoric morphism between central leaves is the composition of a flat universal homeomorphism of finite type and a finite Ã©tale morphism.
\end{Corollary}

\begin{Proof}
  This follows from what has been established about the morphism by using \cite[Lemma~4.8]{messing}.
\end{Proof}

\nocite{mantovan,kottwitz,oort-dim,oort-zink,mantovan-unitary,travaux,haines,zinka,zinkformal,rapo-guide,landvogt-buch,wortmann}

\printbibliography[heading=bibintoc]

\end{document}